\documentclass[a4paper,11pt]{amsart}

\usepackage{epsfig}
\usepackage{color}
\usepackage{amssymb}

\newcounter{lixo}
\newcounter{maislixo}

\newtheorem{theorem}{Theorem}

\newtheorem{lemma}[theorem]{Lemma}
 
\newtheorem{remark}{Remark}  

\newtheorem{Lth}{Theorem} 

\newtheorem{Lcor}[Lth]{Corollary}
\newtheorem{Lprop}[Lth]{Proposition}

\newcommand{\NN}{{\mathbf N}}
\newcommand{\ZZ}{{\mathbf Z}}
\newcommand{\RR}{{\mathbf R}}
\newcommand{\EU}{{\mathbf S}}
\newcommand{\vv}{{\mathbf v}}
\newcommand{\ww}{{\mathbf w}}

\newcommand{\Fc}{{\mathcal F}}
\newcommand{\MM}{{\mathcal M}}
\newcommand{\Vv}{{\mathcal V}}

\newcommand{\Rc}{{\mathcal R}}

\newcommand{\Pp}{{\mathcal P}}

\DeclareMathOperator{\Fix}{Fix}

\newcommand{\dpt}{\displaystyle}

\title[Forcing of a heteroclinic network  \today]{Periodic forcing of a heteroclinic network\\ \today}

\author[I.S. Labouriau and A.A.P. Rodrigues, \today]{
Isabel S. Labouriau
\and Alexandre A. P. Rodrigues}
\address{Centro de Matem\'atica
da Universidade do Porto and\\
Faculdade de
Ci\^encias, Universidade do Porto \\
Rua do Campo Alegre,
687, 4169-007 Porto, Portugal }
\thanks{Centro de Matem\'atica da Universidade do Porto (CMUP --- UID/MAT/00144/2013) is funded by FCT (Portugal) with national (MEC) and European structural funds through the programs FEDER, under the partnership agreement PT2020.   A.A.P. Rodrigues  aknowledges  financial support from Program INVESTIGADOR FCT (IF/00107/2015). Part of this work has been written during AR stay in Nizhny Novgorod University partially supported by the grant RNF 14-41-00044.}
\email{ islabour@fc.up.pt \quad alexandre.rodrigues@fc.up.pt }
\begin{document}

\begin{abstract}
We present a comprehensive mechanism for the emergence of rotational horseshoes and strange attractors in a class of two-parameter families of periodically-perturbed differential equations  defining a flow on a three-dimensional manifold. When both parameters are zero, its flow exhibits an attracting heteroclinic network associated to two periodic solutions. After slightly increasing both parameters, while keeping a two-dimensional connection unaltered, we focus our attention in the case where the two-dimensional invariant manifolds of the periodic solutions do not intersect.
We prove a wide range of dynamical behaviour, ranging from an attracting quasi-periodic torus to rotational horseshoes and H\'enon-like strange attractors. We illustrate our results with an explicit example.

\end{abstract}

    \maketitle


\textbf{2010 --- AMS Subject Classifications} 
{Primary: 37C60;   Secondary: 34C37, 34D20, 37C27, 39A23, 34C28}

\section{Introduction}
A strange attractor   is an invariant set with at least one positive Lyapunov exponent whose basin of attraction has non-empty interior. Nowadays, at least for families of dissipative systems, chaotic dynamics is mostly understood as the persistence of strange attractors
--- occurring for parameters in a set of  positive Lebesgue measure.
 Persistence of dynamics is physically relevant because it means that the phenomenon is \emph{observable} with positive probability. 
Proof of the existence of strange attractors is usually obtained by comparing the dynamics in an invariant set to either the Lorenz or the H\'enon attractors, or through the unfolding of a singularity,
 for a discussion see \cite{BIRR}.

There are few examples of periodically-forced vector fields exhibiting complex dynamics that may be proven analytically. In this article, we give an explicit mechanism to obtain strange attractors in a two-parameter family of vector fields unfolding an attracting heteroclinic network. When the first parameter is different from zero, two normally hyperbolic attracting tori arise near the network. Fixing this parameter and varying the second, the tori break and suspended horseshoes emerge, via the \emph{torus-breakdown} phenomenon \cite{AS91, Anishchenko}. In the meantime, persistent attractors (of H\'enon-type) associated to homoclinic tangencies are created.

Our  \emph{route to chaos}   from an attracting network is different from the routes described by  \cite{Kaneko} and \cite{Bakri}, where the authors use coupled oscillators. Another different itinerary has been described by \cite{FJ2020} in the context of the \emph{Langford system}. These works are discussed in Section \ref{discussion}.

The theory developed in this paper is explicitly applicable to the analysis of various specific differential equations and the results obtained are beyond the capacity of the classical Birkhoff-Melnikov-Smale method associated to \emph{heteroclinic tangles} \cite{Wang_2013}.

Our purpose in writing this paper is not only to point out the range of phenomena that can occur when simple non-linear equations are periodically forced, but to bring to the foreground the techniques that have allowed us to reach these conclusions in a relatively straightforward manner.

We analyse a family of periodic perturbations of an attracting symmetric heteroclinic network defined on the two-sphere. Instead of looking at the time-$T$ maps as in \cite{LR18, TD1}, we extend the phase space and we explicitly compute return maps induced by the perturbed equations in
a neighborhood of the extended heteroclinic network. The  cyclic variable's speed plays an important role in the emergence of the horseshoes.
 Using the techniques explored in  \cite{Shilnikov_book_1, TS86}, we reduce the analysis of the non-autonomous system to that of a two-dimensional map on a circloid.  These techniques are clearly not limited to the systems considered here. It is our hope that they will find applications in other dynamical systems, particularly those that arise naturally from mechanics or physics \cite{AHL2001, DT3, Ruelle, TD1, Rabinovich06}. See also the dynamical description of the periodically-forced van der Pol oscillator in   \cite{Shilnikov_book_1, TS86}, where the authors used the \emph{Afraimovich Annulus Principle} to prove the existence of an invariant torus. 
 
The title \emph{``Periodic forcing of a heteroclinic network''} refers to the study of the dynamics associated to a parametric periodically-forced vector field unfolding  an asymptotically stable heteroclinic network.

\subsection*{Structure of the article}
This article is organized as follows. In Section \ref{sec-setting}, we describe the setting of our problem
and state the main results after having revised some conceptual preliminaries  in Section \ref{preliminaries}. 
Results on general families of vector fields are proved in Section~\ref{secProva_G}, after the derivation  in Section \ref{secPrelim} of the first return map to a given cross section.
 An explicit example is treated in Section~\ref{PropA}.
We finish this article with a discussion on the way our results fit in literature in Section \ref{discussion}.

\section{Preliminaries}
\label{preliminaries}
In this section, we introduce some terminology for vector fields acting on three-dimensional Riemannian manifolds that we will use in the remaining sections.
Consider the two-parameter family of $C^3$--smooth autonomous differential equations
\begin{equation}
\label{general2a}
\dot{x}=F_{(\nu, \mu)}(x)\qquad x\in \EU^3\subset \RR^4  \qquad \nu, \mu\in\RR
\end{equation}
where $\EU^3$ denotes the unit sphere, endowed with the usual topology. Denote by $\varphi_{(\nu, \mu)}(t,x)$, $t \in \RR$, the associated flow. The flow is \emph{complete} (\emph{i.e.}  solutions are defined for all $t\in \RR$) because  $\EU^3$ is a compact manifold without boundary.
 
 \subsection{Heteroclinic structures and symmetry}
Suppose that $A$ and $B$ are two hyperbolic saddles of \eqref{general2a}. There is a {\em heteroclinic cycle} associated to $A$ and $B$ if
$$W^{u}(A)\cap W^{s}(B)\neq \emptyset \qquad \text{and} \qquad W^{u}(B)\cap W^{s}(A)\neq \emptyset.$$ 

The non-empty intersection of $W^{u}(A)$ with $W^{s}(B)$ is called a \emph{heteroclinic connection} between $A$ and $B$, and will be denoted by $[A \rightarrow  B]$. Although the existence of heteroclinic cycles may be a non-generic feature within differential equations, they may be structurally stable within families of systems which are equivariant under the action of a compact Lie group $\mathcal{G}\subset \mathbb{O}(4)$, due to the existence of flow-invariant subspaces (see  \cite{GH}).

Given a group $\mathcal{G}$ of endomorphisms of $\EU^3\subset \RR^4$, we will consider two-parameter families of vector fields $\left(F_{(\nu, \mu)}\right)$ under the equivariance assumption $$F_{(\nu, \mu )}(\gamma x)=\gamma F_{(\nu, \mu )}(x)$$ for all $x \in \EU^3$, $\gamma \in \mathcal{G}$ and $(\nu, \mu )\in  \RR^2.$
For an isotropy subgroup $\widetilde{\mathcal{G}}< \mathcal{G}$, we will write $\Fix(\widetilde{\mathcal{G}})$ for the vector subspace of points that are fixed by the elements of $\widetilde{\mathcal{G}}$. For $\mathcal{G}-$equivariant differential equations, the subspace $\Fix(\widetilde{\mathcal{G}})$ is flow-invariant.

If $\Omega\subset \EU^3$ is a flow-invariant set of \eqref{general2a}, its basin of attraction, $\mathcal{B}(\Omega)$, is the set of points in $\EU^3$ whose orbits have $\omega-$limit in $\Omega$. We say that $\Omega$ is \emph{asymptotically stable} if $\mathcal{B}(\Omega) $ contains 
 all the half-trajectories in positive time starting in
an open neighbourhood of $\Omega$. 
 \subsection{Rotational horseshoes}
Let $\mathcal{H} $ stand for the infinite annulus $\mathcal{H} = \EU^1 \times \RR$ endowed with the usual inner product from $\RR^2$. We denote by $Homeo^+(\mathcal{H} )$ the set of homeomorphisms of the annulus which preserve orientation.
Given a homeomorphism $f :X \rightarrow X$  and a partition of $m\geq 2$ elements $R_1,..., R_{m}$ of $X\subset \mathcal{H}$, the itinerary  function $\xi: X \rightarrow \{1, ..., m\}^\ZZ= \Sigma_m$ is defined by $$\xi(x)(j)=k\quad   \Leftrightarrow \quad f^j(x)\in R_k, \quad \text{for every} \quad j\in \ZZ.$$  
Following \cite{Passeggi}, we say that a compact invariant set $\Lambda \subset \mathcal{H} $ of $f \in Homeo^+(\mathcal{H} )$ is a \emph{rotational horseshoe} if it admits a finite partition $P =\{R_1, ..., R_{m} \}$  by sets  $R_i$ with non empty interior in $\Lambda$ so that
\begin{itemize}
\item 
the itinerary $\xi$ defines a semi-conjugacy between $f|_\Lambda$ and the full-shift $\sigma: \Sigma_m \rightarrow \Sigma_m$, that is $\xi  \circ f = \sigma \circ \xi$ with $\xi$ continuous and onto;
\item 
for any lift $G: \RR^2 \rightarrow \RR^2$ of $f$, there exist  $k>0$ and $m$ vectors $v_1, ...,v_{m} \in \ZZ \times \{0\}$ so that
$$
\left\| (G^n(\hat{x})-\hat{x})  - \sum_{i=0}^n v_{\xi(x)(i)}\right\| <k \qquad \text{for every} \qquad  \hat{x}\in \pi^{-1}(\Lambda), \quad n\in \NN,
$$
where $\|\star\|$ denotes the usual norm of $\RR^2$,  $\pi:\RR^2\rightarrow \mathcal{H}$ denotes the usual projection map and $\hat{x} \in \pi^{-1}(\Lambda)$ is the lift of $x$; more details in the proof of Lemma 3.1 of \cite{Passeggi}. The existence of a rotational horseshoe for a map implies positive topological entropy  at least  $ \log m $.
\end{itemize}

\subsection{Strange attractors}

Strange attractors contribute to the richness and complexity of a dynamical system.  We introduce the following notion, adapted from \cite{MV93}, to the situation under consideration.
A  (H\'enon-type) \emph{strange attractor} of a two-dimensional dissipative diffeomorphism $f$ defined in a Riemannian manifold, is a compact invariant set $\Lambda$ with the following properties:
\begin{itemize}
\item 
$\Lambda$  equals the closure of the unstable manifold of a hyperbolic periodic point;
\item 
the basin of attraction of $\Lambda$   contains an open set (and thus has positive Lebesgue measure);
\item 
there is a dense orbit in $\Lambda$ with a positive Lyapounov exponent (exponential growth of the derivative along the orbit).
\end{itemize}
A vector field possesses a strange attractor if the first return map to a cross section does.

\section{Setting and main results}
\label{sec-setting}
Our object of study is the dynamics around an attracting  heteroclinic network for which we give a rigorous description here.  \subsection{The object of study}
Consider the  family of $C^3$--autonomous differential equations on a two dimensional manifold $\MM^2$ diffeomorphic to the two-sphere $\EU^2\subset \RR^3$ and parametrised by $\nu\in\RR$
\begin{equation}
\label{general0}
 \begin{array}{ll}
\dot x=\Fc_\nu(x), \quad x\in \RR^3.
\end{array} 
 \end{equation}

\begin{figure}
\begin{center}
\includegraphics[height=3.5cm]{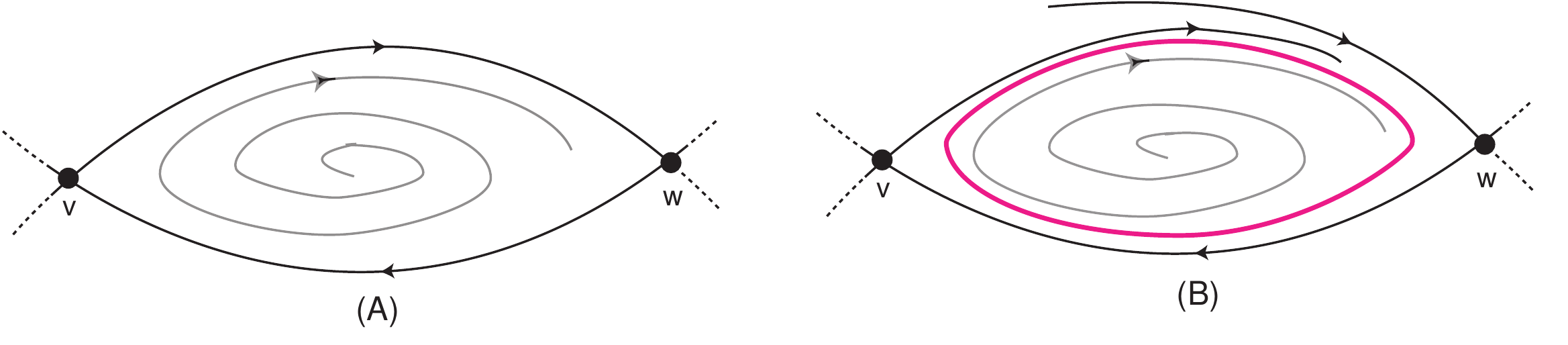}
\end{center}
\caption{\small 
Partial phase portrait associated to  the one-parameter family \eqref{general0},
a periodic solution appears when a heteroclinic connection is broken. (A): $\nu=0$; (B): $\nu\neq 0$.  }
 \label{scheme1}
\end{figure}

Suppose that, for $\nu=0$, the flow of \eqref{general0} has an attracting heteroclinic cycle associated to two 
equilibria
and that for $\nu\ne 0$ one of the heteroclinic connections is broken, yielding an attracting periodic solution, as illustrated in Figure \ref{scheme1}.
More precisely we are assuming that \eqref{general0}   satisfies: 
 
\begin{enumerate}
\renewcommand{\theenumi}{\textbf{(A\arabic{enumi})}}
\renewcommand{\labelenumi}{\textbf{\theenumi}}
\item\label{A1}
  there is  a flow-invariant manifold $\MM^2$ diffeomorphic to $\EU^2$ that is attracting in the sense that every nearby trajectory is asymptotic to it in forward time; 
\item \label{A2}
there are  two equilibria of saddle type
$\vv$ and $\ww$;
\item \label{A3}
there are  two heteroclinic connections from $\ww$ to $\vv$, forming a Jordan curve $\vartheta$ on $\MM^2$.
\setcounter{lixo}{\value{enumi}}
\end{enumerate}
Moreover, the restriction of the flow of \eqref{general0}  to $\MM^2$ satisfies: 
\begin{enumerate}
\renewcommand{\theenumi}{\textbf{(A\arabic{enumi})}}
\renewcommand{\labelenumi}{{\theenumi}}
 \setcounter{enumi}{\value{lixo}}
\item \label{A4}
for $\nu=0$ there are two heteroclinic connections from $\vv$ to $\ww$;
\item \label{A5}
for $\nu=0$ the heteroclinic network formed by the connections between $\vv$ and $\ww$ is attracting;
\item \label{A6}
for $\nu=0$ the only periodic solutions are the equilibria; 
\item\label{A7} for $\nu\ne 0$ the heteroclinic connections from $\vv$ to $\ww$ are broken and there are two attracting hyperbolic periodic solutions, each one in one connected component of   $\MM^2\backslash \vartheta$. 
\setcounter{lixo}{\value{enumi}}
\end{enumerate}

The period of the  solutions of  \ref{A7}  tends to infinity when $\nu$ tends to zero, as they accumulate on the heteroclinic cycles containing $\vv$ and $\ww$. From \ref{A6}, these are the only non-constant periodic solutions. Property~\ref{A7} is the generic unfolding of a heteroclinic network like the one in \ref{A5}. 
In contrast, Property ~\ref{A3} is not generic but occurs generically in systems with a flow-invariant two-dimensional manifold, for instance in the presence of symmetry --- see the example in \S \ref{sec_example} below.

We subject the autonomous differential equation \eqref{general0}  to a family of non-constant  time-periodic perturbations $\phi_\mu(t,x)$ of period $\pi/\omega\in\RR^+$, where the parameter $\mu\ge 0$  controls the amplitude of the perturbation. 
The perturbed equation 
$$
\dot x=\Fc_\nu(x)+\phi_\mu(t,x),\qquad  x\in  \EU^2, 
t\in \RR
$$ 
may be converted into an autonomous equation 
$$
(\dot x,\dot\theta)= F_{(\nu,\mu)} (x,\theta)
$$ 
in $\MM^2\times \EU^1$ by rewriting it as
\begin{equation}
\label{general}
\left\{\begin{array}{ll}
\dot x=\Fc_\nu(x)+\phi_\mu(\theta,x)&{}\\
\dot\theta={2\omega}&\pmod{2\pi}.
\end{array}\right.
\qquad\qquad \phi_0(\theta,x)\equiv 0. 
\end{equation}
 For $\mu=0$ the last equation of \eqref{general} is not coupled to the first. Hence, the dynamics of 
\eqref{general} may be obtained from conditions \ref{A1}--\ref{A7}, as follows.

\begin{Lprop}
\label{PropB1}
If \ref{A1}--\ref{A7} hold for  $\Fc_\nu$
then the flow of \eqref{general} for $\mu=0$ and  small $\nu\ge 0$ satisfies:
\begin{enumerate}
 \item  there is  an invariant flow-invariant manifold $\MM^3$  diffeomorphic to $\EU^2\times \EU^1$
that is globally attracting, in the sense that it attracts all trajectories in its neighbourhood.
\setcounter{maislixo}{\value{enumi}}
\end{enumerate}
Furthermore, the restriction of the flow to $\MM^3$ satisfies:
\begin{enumerate}
 \setcounter{enumi}{\value{maislixo}}
 \item
there are two periodic solutions in $\Pp_\vv=\{\vv\}\times\EU^1$ and  $\Pp_\ww=\{\ww\}\times\EU^1$ of saddle type;
\item
for  $\nu=0$ the invariant manifolds of $\Pp_\vv$ and $\Pp_\ww$ in $\MM^3$ coincide, forming  a  heteroclinic network $\Gamma$ with the geometry of  a  singular two-dimensional torus of genus 2;
\item 
  for $\nu>0$  the manifolds $W^u(\Pp_\ww)$ and $W^s(\Pp_\vv)$  in $\MM^3$ coincide;
\item \label{L5}
  for $\nu>0$ we have  $W^u(\Pp_\vv)\cap W^s(\Pp_\ww)=\varnothing$;
   \item\label{L6}
  for $\nu>0$  there are two attracting  invariant two-dimensional tori $\mathcal{T}^\pm(\nu)$;
 \item\label{L7}
 when  $\nu \rightarrow 0$, the tori $\mathcal{T}^\pm(\nu)$ accumulate on $\Gamma$.
\end{enumerate}
\end{Lprop}

The proof of Proposition \ref{PropB1} follows 
by combining the dynamics of \eqref{general0} with the existence of a cyclic variable (see Chapter 4 of \cite{Shilnikov_book_1}).

\subsection{The periodically forced system}
\label{munot0}

For $\nu=\mu=0$, let $\Sigma$ be a cross section of the heteroclinic cycle $\Gamma$ in $\MM^3$. 
Then $\Sigma$ is also a cross section of \eqref{general} for small $\nu,\mu\ge 0$.
Let $\Rc_{(\nu,\mu)}$ be the first return map to $\Sigma$, with respect to the flow defined by ${F}_{(\nu, \mu)}$.  Define also
$$
\Omega_{(\nu,\mu)}=\left\{X \in \Sigma: \Rc^n_{(\nu,\mu)} (X)\in \Sigma, \quad \forall n \in \NN \right\}\qquad \text{and} \qquad \Lambda_{(\nu,\mu)}= \bigcap_{n \in \NN} \Rc^n_{(\nu,\mu)} \left(\Omega_{(\nu,\mu)}\right).
$$
 In this article, we 
 present a comprehensive analysis on the dynamics 
of $\Rc_{(\nu, \mu)}$ on the \emph{non-wandering set}
 $\Lambda_{(\nu, \mu)}$. When there is no risk of misunderstanding, we omit the subscripts ${(\nu,\mu)}$.
 
When $(\nu,\mu)\ne(0,0)$  one expects that, generically, $W^u(\Pp_\vv) \pitchfork W^s(\Pp_\ww)$.
The case when $W^u(\Pp_\vv)\cap W^s(\Pp_\ww)\ne\varnothing$  has been discussed in \cite{LR17, Wang_2013}, here we are mostly concerned with the case $W^u(\Pp_\vv) \cap W^s(\Pp_\ww)=\varnothing$. 
We also suppose property \ref{A3} (extended to equation \eqref{general}) still holds for the forced system, hence our assumptions are: 
\begin{enumerate}
\renewcommand{\theenumi}{\textbf{(A\arabic{enumi})}}
\renewcommand{\labelenumi}{{\theenumi}}
 \setcounter{enumi}{\value{lixo}}
\item \label{A8}
$W^u(\Pp_\ww) = W^s(\Pp_\vv)$; 
\item \label{A9}
$W^u(\Pp_\vv) \cap W^s(\Pp_\ww)=\varnothing$. 
\end{enumerate}

Our first main result is about the existence  of an invariant   set whose dynamics is conjugate to a full shift over a finite number of symbols. In addition, we also prove the existence of observable chaos.

 \begin{Lth} 
  \label{role_omega}
  If \ref{A1}--\ref{A9} hold for \eqref{general} then
  \begin{enumerate}
  \item for every small $\nu,\mu>0$ there exists $\omega_0>0$ such that if $\omega>\omega_0$, then  $\Lambda_{(\nu, \mu)}$ contains  an invariant set whose dynamics is conjugate to a full shift in two symbols;
  \item 
for $(\nu, \mu)$ in a set  $\mathcal{U}\subset \RR^2$  of positive Lebesgue measure, the return map
$\Rc_{(\nu,\mu)}$ exhibits a strange attractor.
\end{enumerate}
 \end{Lth}
 
   The dynamics of  $\Lambda_{(\nu, \mu)}$ is mainly governed by the geometric configuration of the global invariant manifold $W^u(\Pp_\vv)$.
  The proof of  this result is  
  done is Section~\ref{secProva_G}.  Horseshoes of Theorem  \ref{role_omega} have a different nature from those associated to the \emph{heteroclinic tangle} in which the manifolds have a transverse intersection  \cite{LR17, Wang_2013}.  This will be discussed in Section \ref{discussion}.

For a fixed $\nu>0$, if the ratio of $\omega$ and the period of the hyperbolic periodic solution of  \ref{A7}  is irrational, then trajectories on the torus $\mathcal{T}^\pm(\nu)$ are unlocked, in the sense that they never close. These solution on  the torus are called \emph{quasiperiodic} \cite{Herman, Shilnikov_book_1}. 
  If the frequencies have a rational ratio,  trajectories are locked.

In a resonant torus, where all solutions are locked, the frequency locking ratio $p/q$ means that while the $x$ component of  a solution turns $p$ its $\theta$ component winds $q$ times.
 This ratio is related  to the \emph{rotation number} associated to the periodic orbit \cite{Herman} and will be used in the proof of the second part of Theorem \ref{role_omega}.

\subsection{The example}
\label{sec_example}
 An explicit two-parameter family $\Fc_\nu(x)+\phi_\mu(t,x)$ of vector fields   in $\EU^2\subset \RR^3$  such that $\Fc_\nu(x)$ satisfies  \ref{A1}--\ref{A7}  is given by
 \begin{equation}
\label{general4}
\left\{ 
\begin{array}{l}
\dot x_1 = x_1(1-r^2)-\alpha x_1 x_3 +\beta x_1x_3^2 + (1- x_1) \,\, \,  (\mu\,  [f (\theta)-1]+\nu)\\
\dot x_2 =  x_2(1-r^2) + \alpha x_2 x_3 + \beta x_2 x_3^2 \\
\dot x_3 = x_3(1-r^2)-\alpha(x_2^2-x_1^2)-\beta x_3 (x_1^2+x_2^2) \\
 \dot\theta={2\omega}\pmod{2\pi}
\end{array}
\right.
\end{equation}
where 
$$
\nu, \omega\in \RR^+\quad \mu\in\RR\qquad
r^2= x_1^2+x_2^2+x_3^2, \qquad \beta<0<\alpha, \qquad
  |\beta|<\alpha ,
$$
 and $f$ is a  non constant $2\pi$-periodic map of class $C^3$.

 For $\mu=\nu=0$, the equation $\dot x=\Fc_{0}(x)$, $x\in \RR^3$,  is one of the examples constructed and analysed in   \cite{ACL06} and also  studied in \cite{LR18}.  
 The perturbing  term $(1-x_1) \  [\mu (f(2\omega t) -1)+\nu] $ appears only in the first coordinate
  for two reasons. First, it simplifies the computations. Secondly, it allows comparison with previous work by other authors \cite{AHL2001,DT3, Rabinovich06, TD1}. 

  \begin{Lprop}
\label{periodic_solution_prop}
 The vector field $\Fc_\nu$ associated to  \eqref{general4} at $\mu=0$ is equivariant under the action of $\kappa(x,y,z)=(x, -y,z)$ 
and therefore the plane $\Fix(\ZZ_2(\kappa))=(x,0,z)$ is flow-invariant.
If   $|\nu|>0$ is small,  the flow of $\dot\zeta=\Fc_\nu(\zeta)$ satisfies conditions \ref{A1}--\ref{A7}.
In particular, the flow-invariant curve $\MM^2\cap \Fix(\ZZ_2(\kappa))$ consists of  two equilibria of saddle-type $\vv$ and $\ww$ and two heteroclinic connections from $\ww$ to $\vv$. 
There are also four equilibria in $\MM^2$  that are repelling foci and these are all the equilibria in $\MM^2$.
\end{Lprop}

 The proof of  this result is the content of Section~\ref{PropA}.

 From \eqref{L6} in Proposition~\ref{PropB1} it follows that there are two invariant tori for the flow of the equation $(\dot x,\dot\theta)= F_{(\nu,\mu)} (x,\theta)$
 associated to  \eqref{general4}.
 Existence of invariant tori  is usually shown using the Afraimovich Annulus Principle \cite{AS91}, here we show it directly by reducing the problem to a two-dimensional manifold and applying the Poincar\'e-Bendixson theorem.
 
 From Theorem~\ref{role_omega} and Proposition~\ref{periodic_solution_prop} it follows immediately
 
 \begin{Lcor}
 \label{CorExemploCaos}
For small $\mu>0$, $\nu>0$ and  for $\omega>0$ large enough,  the flow of \eqref{general4} exhibits a hyperbolic rotational horseshoe
and for a set of positive Lebesgue measure of parameters it also contains  strange attractors. 
 \end{Lcor}

\section{Local coordinates and first return map}
\label{secPrelim}

In this section we will analyse the dynamics near the  heteroclinic attractor $\Gamma$ through local maps, after selecting appropriate coordinates near the saddles   $\Pp_\vv=\{\vv\}\times\EU^1$ and  $\Pp_\ww=\{\ww\}\times\EU^1$.
\subsection{Geometry near $\Pp_{\vv}$ and $\Pp_{\ww}$}
\label{Local}
Let $U_a$
be pairwise disjoint compact neighbourhoods in $\MM^3$ of the nodes  $\Pp_{a}$, $a\in \{\vv,\ww\}$, such that each boundary  $\partial U_a$ is a finite union of smooth surfaces delimited by curves, each surface
 transverse to the vector field everywhere, except at its
 boundary.  
Each  $U_a$ is called an \emph{isolating block} for  $\Pp_{a}$ and, topologically, it consists of a hollow cylinder.  

\begin{figure}
\begin{center}
\includegraphics[height=9cm]{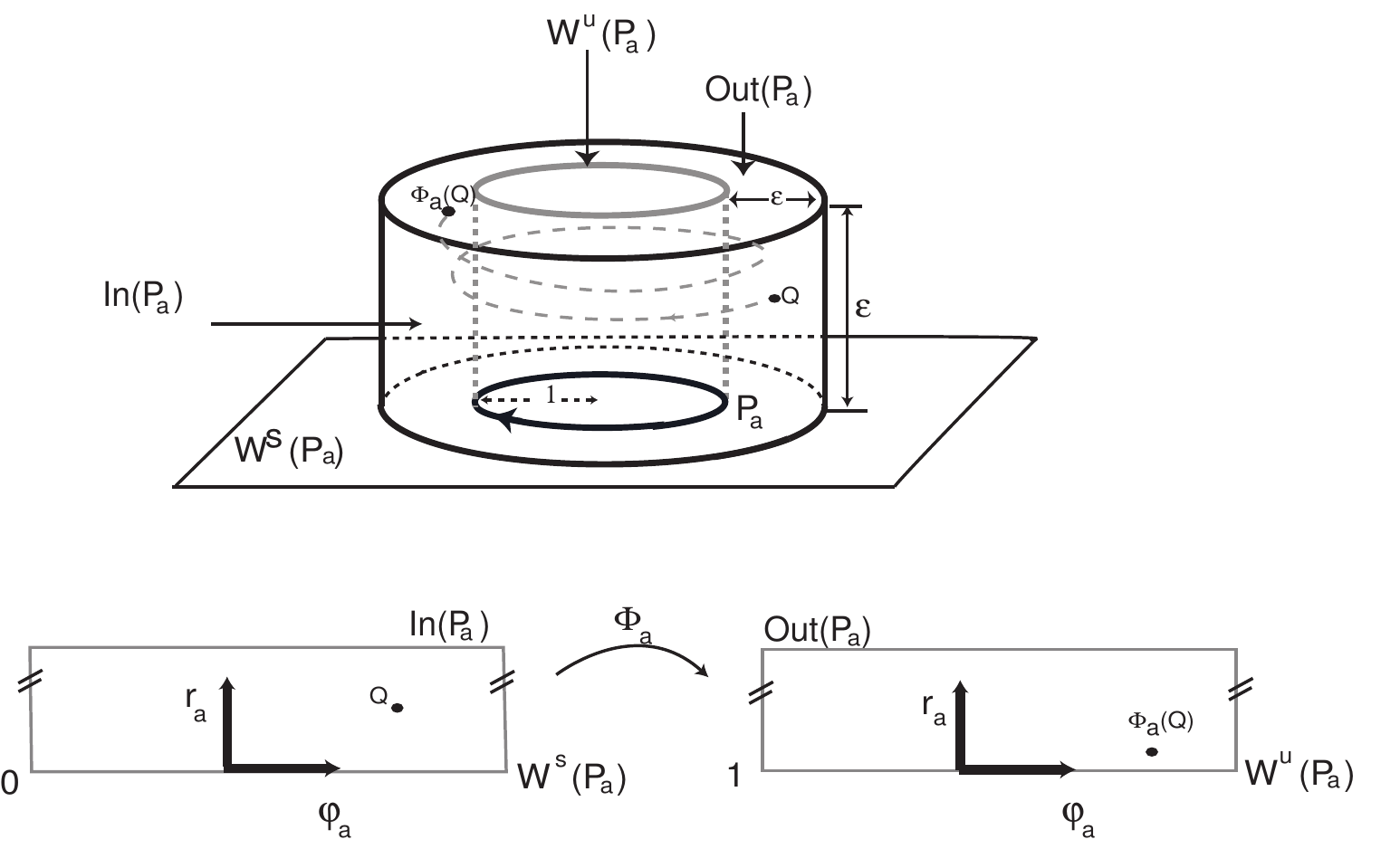}
\end{center}
\caption{\small Top:  isolating block near the periodic solution $\Pp_a$, $a\in \{\vv,\ww\}$.
Bottom: coordinates at the $In$ and $Out$ components of the boundary.  
Double bars mean that the sides are identified.    }
\label{local_map_LR_torus1}
\end{figure}

For $a\in \{\vv,\ww\}$, let $\Sigma_a$ be a cross section to the flow at $p_a \in \Pp_{a}$. Since $\Pp_{a}$ is hyperbolic, there is a neighbourhood $U^*_a$ of $p_a$ in $\Sigma_a$ where the first return map to $\Sigma_a$
is $C^1$ conjugate to its linear part. 
 Let $e^{-c_a}$ and $e^{e_a}$, with $c_a,e_a>0$,  be the eigenvalues
of the derivative $D\Fc_{(\nu, \mu)}(a)$.
Then, for each $k\ge 2$ there is an open and dense subset of $\RR^2$ such that, if 
$(-c_a,e_a)$
 lies
 in this set, then the conjugacy is of class $C^k$ (details may be checked in  Appendix A of \cite{LR17}).

Suspending the linear map gives rise, in cylindrical coordinates $(\rho, \theta, z)$ around $\Pp_{a}$, to the  equations
\begin{equation}
\label{ode of suspension}
\left\{ 
\begin{array}{l}
\dot{\rho}=-c_{a}(\rho -1) \\ 
\dot{\theta}=2\omega \\ 
\dot{z}=e_{a}z
\end{array}
\right.
\end{equation}
 whose flow is
$C^2$-conjugate to the original flow near $\Pp_{a}$. In these coordinates, the periodic trajectory $\Pp_{a}$ is the circle defined by $\rho=1$ and $z=0$.
For the moment, let $W^s_{loc}(\Pp_{a})$ and $W^u_{loc}(\Pp_{a})$ be the connected components of $W^s(\Pp_{a})$ and $W^u(\Pp_{a})$, respectively, contained in the suspension of $U^*_a$  and containing $\Pp_a$ in their closure.
In these coordinates,  $W^s_{loc}(\Pp_{a})$, is  the plane $z=0$ and   $W^u_{loc}(\Pp_{a})$ is the surface $\rho=1$.

As illustrated in Figure \ref{local_map_LR_torus1}, we consider a  hollow three-dimensional cylinder $V_a(\varepsilon_a)$  of 
$\Pp_{a}$ contained in the suspension of $U^*_a$ 
(with small $\varepsilon_a>0$ to be determined later)
given by
$$
V_a(\varepsilon_a)=\left\{ (\rho,\theta,z):\quad 1\le\rho< 1+\varepsilon_a,
\quad 0\le z< \varepsilon_a\quad \text{and}\quad 
\theta\in\RR\pmod{2\pi}
\right\}\  .
$$
When there is no ambiguity, we write $V_a$ instead of $V_a(\varepsilon_a)$.
Its boundary contains the trajectory $\Pp_{a}$ and is a 
union
$$
\partial V_{a}= In(\Pp_{a}) \cup Out(\Pp_{a}) \cup \mathcal{W} (\Pp_{a})
$$
where 
\begin{itemize}
\item
$\mathcal{W} (\Pp_{a})=\Pp_{a}\cup  \left(W^s_{loc}(\Pp_{a})\cap V_a\right) \cup  \left(W^u_{loc}(\Pp_{a})\cap V_a\right)$.
\item
$ W^s_{loc}(\Pp_{a})\cap V_a$ is the lower boundary of the hollow cylinder, given by $z=0$, $1<\rho\le 1+\varepsilon_a$.
\item
$W^u_{loc}(\Pp_{a})\cap V_a$ is the inner boundary of the hollow cylinder, given by $\rho=1$, $0<z\le\varepsilon_a$.
\item 
$In(\Pp_{a})$ is the outer  wall of the cylinder, defined by $\rho=1+\varepsilon_a$, $0\le z\le\varepsilon_a$.\\
Trajectories starting at  $In(\Pp_{a})$ go inside $V_a$ in small positive time.
\item 
$Out(\Pp_{a})$ is the  top  of the cylinder, the annulus defined by $z=\varepsilon_a$, $1\le\rho\le 1+\varepsilon_a$.\\
Trajectories starting at $Out(\Pp_{a})$ go inside $V_a$ in small negative time. 
\item 
The vector field is transverse to $In(\Pp_{a})\cup Out(\Pp_{a})$ except at the circle 
$In(\Pp_{a})\cap Out(\Pp_{a})$.  
\end{itemize}

The cylinder wall $In(\Pp_{a})$ is
parametrised by the covering map
$$
(\varphi,r)\mapsto(1+\varepsilon_a,\varphi,r)=(\rho,\theta,z),
$$
where $\varphi\in\RR$ and
$0\le r\le\varepsilon_a$. 
The annulus $Out(\Pp_{a})$ is parametrised by the covering
$$
(\varphi,r) \mapsto ( r,\varphi, \varepsilon_a)=(\rho,\theta,z),
$$
for $1\le r\le 1+\varepsilon_a$ and $\varphi\in\RR$.

For $\nu\neq 0$,  $\mu= 0$ we have from \eqref{L5} of Proposition~\ref{PropB1} that $W^u(\Pp_\vv) \cap W^s(\Pp_\ww)=\varnothing$ and the same assumption is made for $\mu\ne 0$ in Theorem~\ref{role_omega}.
Since  for $\nu=\mu=0$ these invariant manifolds coincide, then for small $\mu,\nu\ne 0$ the manifold $W^u(\Pp_\vv)$ must come close to $\Pp_\ww$.
Without loss of generality we may assume that it meets $In(\Pp_\vv)$,
we are concerned with the parameter values for which this holds.

From now on the portion of the unstable manifold of $\Pp_\vv$  that goes from $\Pp_\vv$ to $In(\Pp_\ww)$ without intersecting $V_\ww$ will be  denoted $W^u_{loc}(\Pp_\vv)$. Similarly, $W^s_{loc}(\Pp_\ww)$ will denote the portion of the stable manifold of $\Pp_\ww$ that is outside $V_\ww$ and goes directly from $Out(\Pp_\ww)$ to $\Pp_\vv$ in negative time.

\subsection{Local and global maps}\label{sublocal}
\begin{figure}
\begin{center}
\includegraphics[height=5cm]{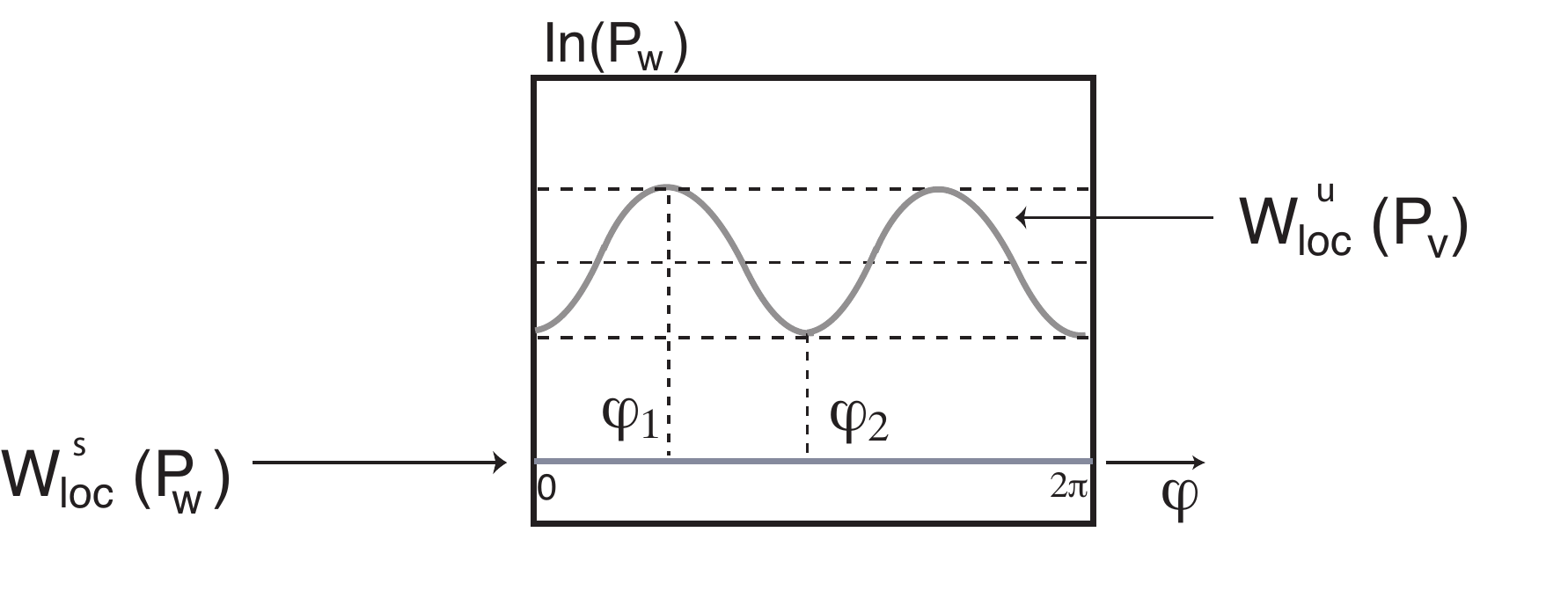}
\end{center}
\caption{\small The local unstable manifold of  $\Pp_{\vv}$ intersects the cylinder $In(\Pp_{\ww})$  on a closed curve, the graph of a periodic Morse function  $\xi_{(\nu, \mu)}:\RR \rightarrow [0, \varepsilon_\ww]$ with $\nu>0$, $\mu>0$. }
\label{xi1n}
\end{figure}

For each $a\in\{\vv, \ww\}$, we may solve  (\ref{ode of suspension}) explicitly, then we compute the flight time from  $In(\Pp_{a})$ to $Out(\Pp_{a})$ by solving the equation $z(t)=\varepsilon_a$ for the trajectory whose initial condition is 
$(\rho,\theta,z)=(1+\varepsilon_a,\varphi,r) \in In(\Pp_{a})\backslash W^s(\Pp_{a})$, with $z>0$, as in \cite{LR17}.
  Replacing this time in the other coordinates of the solution, yields the local map 
  $$
  \Phi _{a}:In(\Pp_{a})\backslash W^s(\Pp_{a})\, \longrightarrow\,  Out(\Pp_a)
  $$
  given by
 \begin{equation}
\label{local map}
\Phi _{a}(\varphi,r)=
\left(\varphi-\frac{2\omega}{e_a}\ln\left(\frac{r}{\varepsilon_a}\right),
1+ \varepsilon_a \left(\frac{r}{\varepsilon_a}\right)^{\delta_a}\right) 
 \end{equation}
 where $\delta_a=\dfrac{c_{a}}{e_{a}}>0$. 
For the transition maps from one isolating block to the other we use assumptions \ref{A8} and \ref{A9}.
With this notation, we formulate them as follows: 
 \begin{itemize}  
\item 
The two sets  $W^u(\Pp_{\ww})$  and $W^s(\Pp_{\vv})$ coincide; 
\item
The manifold $W^u_{loc}(\Pp_{\vv})$ intersects
 the cylinder $In(\Pp_{\ww})$  on a non-contractible closed curve $\gamma_{(\nu, \mu)}$. 
\end{itemize}
We will assume that  the universal cover of $\gamma_{(\nu, \mu)}$ is the graph of a smooth
Morse function  
$$
\xi_{(\nu, \mu)}:\RR \rightarrow [0, \varepsilon_\ww],
$$
 as in Figure \ref{xi1n},
satisfying the following conditions for $\nu>0$: 
\begin{itemize}
\item 
$\xi_{(\nu, \mu)}$ is not constant 
 (because the perturbing term $\phi_\mu$ is not constant)
 and is $2\pi$-periodic; 
\item 
$\xi_{(\nu, \mu)}$ has a local maximum at  $\varphi_1>0$, a local minimum at $\varphi_2>\varphi_1$ and no other critical point in the interval $(\varphi_1,\varphi_2)$;  
\item if $\mu>0$ and $\nu>0$ then $\forall \varphi\in  \RR$,  $\xi_{(\nu, \mu)}(\varphi)>0$;  
\item $\dpt \lim_{\nu\rightarrow 0}\,  \max_{\varphi\in \RR } \xi_{(\nu, \mu)}(\varphi)=0$. 
\end{itemize}

With these assumptions, we may take  the maps $$\Psi_{\vv \rightarrow \ww}: Out(\Pp_\vv)\longrightarrow In (\Pp_\ww) \qquad \text{and} \qquad \Psi_{\ww \rightarrow \vv}: Out(\Pp_\ww)\longrightarrow In (\Pp_\vv)$$ to be given by
\begin{equation}\label{transition21}
\qquad \Psi_{\vv \rightarrow \ww}(\varphi,r)=\left(\, \varphi, \,\, (r-1)+ \xi_{(\nu, \mu)}(\varphi)\right) \qquad\text{and}\qquad \Psi_{\ww \rightarrow \vv}(\varphi,r)=\left(\, \varphi, \,\, (r-1)\right).
\end{equation}
\subsection{The return map}
 Let  $\Rc_{(\nu,  \mu)}=\Phi_{\vv} \circ \Psi_{\ww \rightarrow \vv} \circ \Phi_{\ww} \circ \Psi_{\vv \rightarrow \ww}$
  be the first return map to $Out(\Pp_\vv)$, well defined on the set   of initial conditions $(\varphi,r) \in Out(\Pp_\vv)$ whose solution returns to $Out(\Pp_\vv)$. 
 For  $r>1$, the map $\Rc_{(\nu,\mu)}$ is given by
 \begin{eqnarray*}\label{first1}
\Rc_{(\nu, \mu)}(\varphi,r)&=& \left[ 
\varphi-\omega K\ln \left[(r-1)+\xi_{(\nu, \mu)} (\varphi)\right] -\omega k_\varepsilon\pmod{2\pi},\ 
 1+\dfrac{\varepsilon_\vv}{\varepsilon_\ww^{\delta_\ww}}\left[(r-1)+\xi_{(\nu, \mu)} (\varphi) \right] ^\delta\right]\\
&=& \left(R_1(\varphi,r), R_2(\varphi ,r)\right)
\end{eqnarray*}
where
$$
\delta = \delta_\vv \delta_\ww>1, \qquad  
K = 2 \left( \frac{e_\vv + c_\ww}{e_\vv \, e_\ww} \right)>0 \qquad \text{and}\qquad 
k_\varepsilon =-2K\ln \varepsilon_\ww>0.
$$
The map $\Rc_{(\nu, \mu)}$ is well defined if $\varepsilon_\vv+\dpt\max_{0\le \varphi\le 2\pi}\xi_{(\nu, \mu)} (\varphi)\le \varepsilon_\ww$.
In particular, we need $\varepsilon_\vv< \varepsilon_\ww$. 
 The inequality $\delta>1$ comes from Property  \ref{A5} and the conditions  of \cite{KM1, KM2} for a heteroclinic cycle to be attracting.

\section{Proof of Theorem~\ref{role_omega}}
\label{secProva_G}

The goal of the section
is to obtain  an invariant set  $\Lambda\subset Out^\pm(\Pp_\vv)$ where the map $\Rc_{(\nu, \mu)}|_\Lambda$ is topologically conjugate to a Bernoulli shift with two symbols. 
The argument uses the Conley-Moser conditions, see for instance \cite{Wiggins}.
\subsection{Stretching the angular component}
Let $[\varphi_L,\varphi_R]$ be an interval where $\xi_{(\nu, \mu)} (\varphi)$ is monotonically decreasing, with 
$$
\xi_L=\xi_{(\nu, \mu)} (\varphi_L)>\xi_{(\nu, \mu)} (\varphi_R)= \xi_R
$$
 and consider $\mathcal{D}\subset Out(\Pp_{\vv})$ parametrised by $(\varphi,r)\in [\varphi_L,\varphi_R]\times[1,1+\varepsilon_\vv]$, with $\varepsilon_\vv+\xi_L<\varepsilon_\ww$.
Then $\mathcal{D}$ is a set  of initial conditions $(\varphi,r) \in Out(\Pp_\vv)$ whose solution returns to $Out(\Pp_\vv)$. We start by establishing some properties of the map $\Rc $  within this set.

\begin{lemma}
\label{lema:contract}
For small $\nu,\mu>0$, the following assertions hold in  $\mathcal{D}$  with $\varepsilon_\vv+\xi_L<\varepsilon_\ww$: 
\begin{enumerate}
\item
for any $r\in [1,1+\varepsilon_\vv]$ the map $\varphi\to R_1(\varphi,r)$ is an expansion;
\item
for any $\varphi\in[\varphi_L,\varphi_R]$  the  map $r\to R_2(\varphi,r)$ is a contraction. 
\end{enumerate}
\end{lemma}

\begin{proof} 
The first assertion follows from
$$
\dfrac{\partial  R_1 (\varphi,r)}{\partial \varphi } =
1-\dfrac{\omega K}{r-1+\xi_{(\nu, \mu)}( \varphi) }\dfrac{d\xi_{(\nu, \mu)}( \varphi)}{d\varphi}>1
$$
because $\xi_{(\nu, \mu)}( \varphi) >0$ and  $\dfrac{d\xi_{(\nu, \mu)}( \varphi)}{d\varphi}<0$ since we are assuming $\xi_{(\nu, \mu)}$ is monotonically decreasing in $[\varphi_L,\varphi_R]$.
The second assertion follows from
$$
\dfrac{\partial  R_2 (\varphi,r)}{\partial r } = \delta\dfrac{\varepsilon_\vv}{\varepsilon_\ww^{\delta_\ww}} \left[(r-1) +\xi_{(\nu, \mu)}( \varphi)\right]^{\delta-1}. 
$$
Since  $0<\xi_{(\nu, \mu)}( \varphi)\le \xi_L$  and $0<r-1\leq \varepsilon_\vv\le\varepsilon_\ww<1$ then
$$
0<\frac{\partial  R_2 (\varphi,r)}{\partial r } 
\le \delta\dfrac{\varepsilon_\vv}{\varepsilon_\ww^{\delta_\ww}}\left[\varepsilon_\ww+\xi_L \right]^{\delta-1}
=\delta\dfrac{\varepsilon_\vv}{\varepsilon_\ww}+\mathcal{O}(\xi_L),
$$
where $\mathcal{O}$ stands for the standard Landau notation. We have $\dpt\lim_{\nu\rightarrow 0} \max_{\varphi\in[\varphi_L,\varphi_R]} \xi_{(\nu, \mu)}(\varphi)=0$,
therefore $0<\dfrac{\partial  R_2}{\partial r} (\varphi,r) <1 $ for small $\mu>0$. 
\end{proof}

We call the graph $(\varphi,s(\varphi))$ in  $\mathcal{D}$ of a monotonic map $s(\varphi)$ with 
$\varphi_L\le\varphi\le\varphi_R$, a  \emph{segment across} $\mathcal{D}$.

\begin{figure}
\begin{center}
\includegraphics[width=4cm]{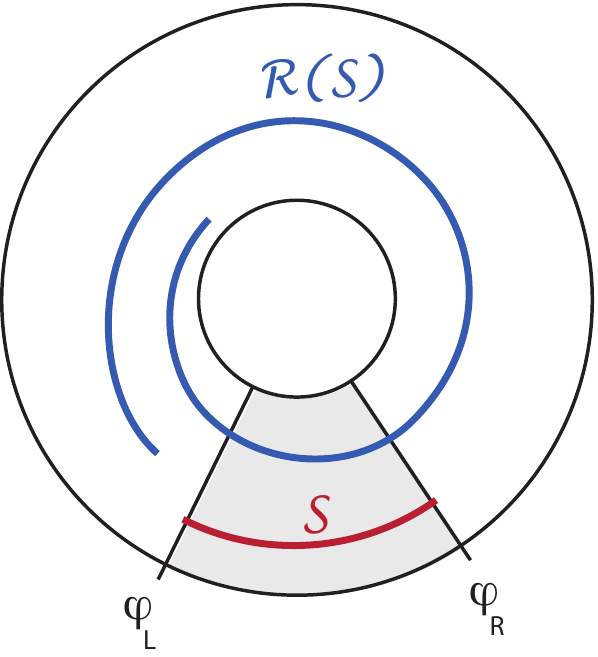} 
\end{center}
\caption{\small When $\omega\ge\omega_0$, the segment $S$ (red) in the domain  $\mathcal{D}\subset Out(\Pp_{\vv})$ (gray) is  transformed by the first return map $\Rc $ into a curve (blue) that makes a full turn around $Out(\mathcal{P}_\vv)$ intersecting $\mathcal{D}$ in at least one segment.}
\label{fig:spiral}
\end{figure}

\begin{lemma}\label{lema:spiral}
Consider the segment 
$S=\{(\varphi,r_*)\ \varphi_L\le\varphi\le\varphi_R \}\subset \mathcal{D}$
for a given $r_*>0$.
For  small $\mu,\nu>0$   with $\varepsilon_\vv+\xi_L<\varepsilon_\ww$, if 
$$
\omega\ge\omega_0=\dfrac{ 2\pi}{K\ln( 1+(\xi_L-\xi_R)/(1+\xi_R))}
$$ 
then for any $r_*\in(1,1+\varepsilon_\vv]$,
the set $\Rc (S)\cap \mathcal{D}$ is a curve   containing a segment across $\mathcal{D}$.
\end{lemma}

\begin{figure}
\begin{center}
\includegraphics[height=6cm]{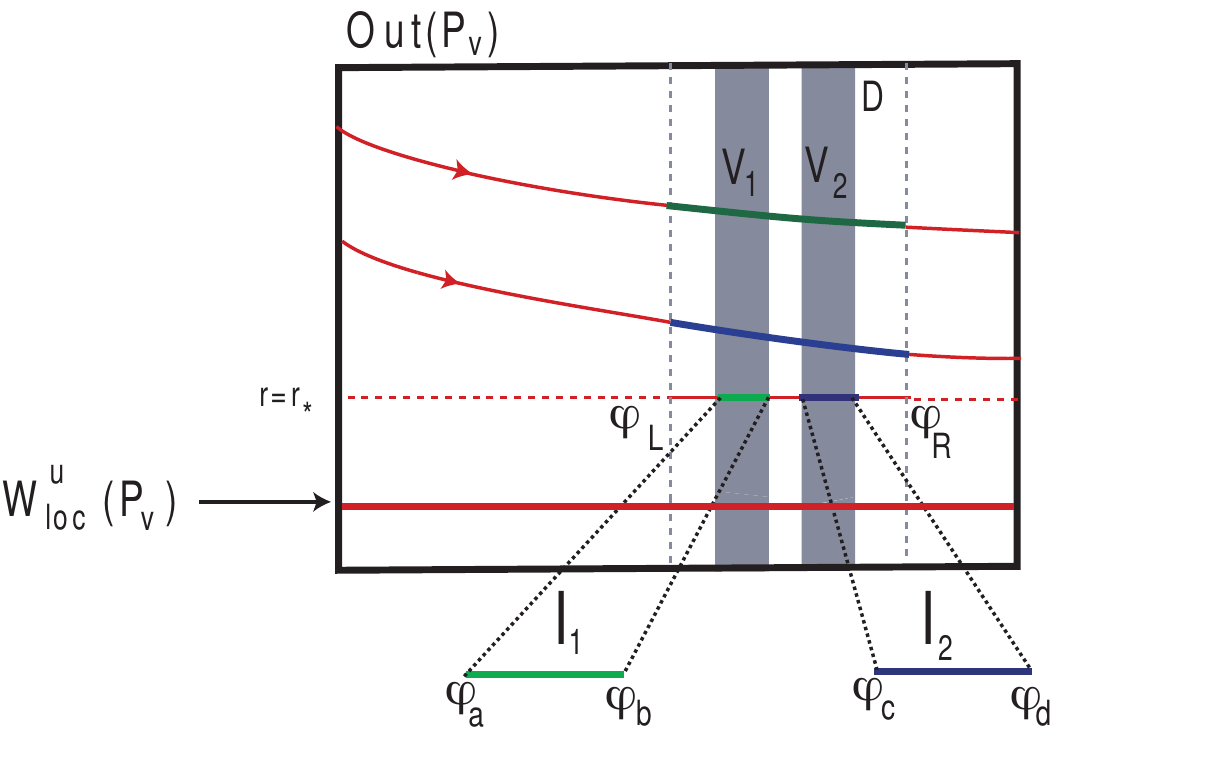}
\end{center}
\caption{\small  For $r_\star \in  [0, 1]$, the image under $\Rc$ of the segment $[\varphi_L, \varphi_R]\times \{r=r_* \}$ is a curve (without folds) intersecting twice the rectangle $\mathcal{D}$. 
The  subsets $I_1$ and $I_2$ are stretched by the first return map $\Rc$ into segments across $\mathcal{D}$. }\label{horseshoe_fig}
\end{figure}

\begin{proof}
Since  $\xi_{(\nu, \mu)} (\varphi)$ is monotonically decreasing in $[\varphi_L,\varphi_R]$ then  the map 
$ \varphi\to R_2(\varphi, r_*)$ is monotonically decreasing and 
$$ \varphi\mapsto \varphi-\omega K\ln \left[(r-1)+\xi_{(\nu, \mu)} (\varphi)\right] -\omega k_\varepsilon=R_1(\varphi,r_*)$$ is monotonically increasing  in the same interval. 
Hence $\Rc (S)$ is the graph of a monotonic map $s(\varphi)$, with the map $s$ defined in some interval  $I$, as in Figure~\ref{fig:spiral}. 
It remains to obtain an estimate of the variation of the first coordinate of $\Rc (S)$ to ensure that $[\varphi_L,\varphi_R]\subset I$.
From the definition of $R_1$ and properties of the logarithm, one knows that the difference $\Delta =R_1(\varphi_R,r_*) - R_1( \varphi_L, r_*)$ satisfies
$$
\begin{array}{rl}
\Delta =&
\left(\varphi_R-\omega K\ln \left[(r_*-1)+\xi_{(\nu, \mu)} (\varphi_R)\right] -\omega k_\varepsilon\right)-
\left(\varphi_L-\omega K\ln \left[(r_*-1)+\xi_{(\nu, \mu)} (\varphi_L)\right] -\omega k_\varepsilon\right)\\ \\
=&(\varphi_R-\varphi_L)+\omega K
 \ln\dfrac{(r_*-1)+\xi_{(\nu, \mu)} (\varphi_L)}{(r_*-1)+\xi_{(\nu, \mu)} (\varphi_R)}\\ \\
 =&(\varphi_R-\varphi_L)+\omega K
  \ln \dfrac{(r_*-1)+\xi_L}{(r_*-1)+\xi_R}> (\varphi_R-\varphi_L)
\end{array}
$$
where for the last inequality we use $\xi_L>\xi_R$ hence $ \ln \dfrac{(r_*-1)+\xi_L}{(r_*-1)+\xi_R}>0$.
Moreover,
$$
\dfrac{(r_*-1)+\xi_L}{(r_*-1)+\xi_R}=
1+\dfrac{\xi_L-\xi_R}{(r_*-1)+\xi_R}\ge 
1+\dfrac{\xi_L-\xi_R}{1+\xi_R}.
$$
Therefore, if $\omega\ge\dfrac{ 2\pi}{K\ln( 1+(\xi_L-\xi_R)/(1+\xi_R))}$, then $\Delta \ge 2\pi+  (\varphi_R-\varphi_L)$ and hence the curve $\Rc (S)$ goes across $\mathcal{D}$ at least once, as in Figures~\ref{fig:spiral} and \ref{horseshoe_fig}.
\end{proof}

\subsection{Proof of Theorem~\ref{role_omega}. Part I}
\label{prova_parte1}
Given a rectangular region   in $Out(\Pp_\vv)$, parametrised by 
$ [\varphi_a, \varphi_b]\times [r_1,r_2]$, a \emph{vertical  strip} in the region is a set
$$
\Vv=\{(\varphi,r): \varphi\in[u_1(r),u_2(r)]\qquad r \in \,[r_1,r_2]\}
$$
where $u_1,u_2: [r_1,r_2] \rightarrow [\varphi_a,\varphi_b]$ are Lipschitz functions with Lipschitz constants less than $\mu_v\ge 0$, such that $u_1(r)<u_2(r)$. The \emph{ vertical  boundaries} of a  vertical  strip are the graphs of the maps $u_i$; the \emph{horizontal boundaries} are the lines $\{r_i\} \times  [u_1(r_i),u_2(r_i)]$, $i=1,2$; 
the \emph{width} $d(\Vv)$ of the strip is 
$$
d(\Vv)=\max_{r_1\le r\le r_2} |u_1(r)-u_2(r)|.
$$
In an analogous way we define a \emph{horizontal strip} across  a \emph{horizontal rectangle} in $Out(\Pp_\vv)$ with the roles of $\varphi$ and $r$ reversed, Lipschitz constants less than $\mu_h\ge 0$.

Let $I_1=[\varphi_a, \varphi_b]$ and $I_2=[\varphi_c, \varphi_d]$ be two disjoint intervals satisfying
$$
\varphi_L\le\varphi_a<\varphi_b<\varphi_c<\varphi_d\le \varphi_R.
$$
We claim that for $\omega\ge\omega_0$ (of Lemma~\ref{lema:spiral}) the  two  vertical  strips
$$
\Vv_i= \{(\varphi,r)\in \mathcal{D}: \varphi\in I_i\}\qquad i=1,2
$$
satisfy the Conley-Moser conditions \cite{Wiggins}: 
\begin{enumerate}
\renewcommand{\theenumi}{(P\arabic{enumi})}
\renewcommand{\labelenumi}{{\theenumi}}
\item\label{f(V)=H}
The image $\Rc (V_i)\cap \mathcal{D}=H_i$ is the union of disjoint horizontal strips across $\mathcal{D}$.
\item\label{f(V)capV}
For every vertical strip $V\subset (V_1\cup V_2)$ the set $\Rc ^{-1}(V)\cap V_i=\widetilde{V_i}$ is a vertical strip across $\mathcal{D}$ with $d\left(\widetilde{V_i}\right)\le \lambda_v d(V)$ for some $\lambda_v\in(0,1)$.
\item\label{f(H)capH}
For every horizontal strip $H\subset (H_1\cup H_2)$ the set $\Rc (H)\cap H_i=\widetilde{H_i}$ is a horizontal strip across $\mathcal{D}$ with $d\left(\widetilde{H_i}\right)\le \lambda_h d(H)$ for some $\lambda_h\in(0,1)$.
\end{enumerate}
In order to establish the claim, note that by Lemma~\ref{lema:spiral} each line with constant $r$ in $\mathcal{V}_i$ is mapped by $\Rc$ into a curve in $Out(\mathcal{P}_\vv)$ that intersects $\mathcal{D}$ in at least one segment. 
 Lemma~\ref{lema:contract} ensures that as $r$ varies in $]1,1+\varepsilon_\vv]$, the second coordinate of their image varies in an interval of length less than $\varepsilon_\vv$, so the union of these segments lies in horizontal strips across $\mathcal{D}$, establishing \ref{f(V)=H}.
Properties~\ref{f(V)capV} and  \ref{f(H)capH} follow from Lemma~\ref{lema:contract}.
From this claim it follows that there exists an $\Rc$-invariant set of initial conditions 
$$
{\Lambda}=\bigcap_{n \in \ZZ} \Rc^n(\Vv_1 \cup \Vv_2)
$$ 
on which $\left.\Rc\right|_\Lambda$ is topologically conjugate to a Bernoulli shift on two symbols. By construction it is a rotational horseshoe (according to \cite{Passeggi}) with $m=2$.

\subsection{Proof of Theorem~\ref{role_omega}. Part II}
For $\nu>0$ and $\mu=0$, the flow of \eqref{general} has an attracting two-dimensional torus (by Proposition \ref{PropB1}). In particular, there is a cross section $\Sigma$ where the torus defines an invariant curve $\mathcal{C}$ under the first return map $\Rc$.  Furthermore, there is countable set of values of the type $(\nu_i, 0)$, $i \in \NN$, for which the first return map $\Rc $  has at least one saddle and a sink lying on $\mathcal{C}$ ($\Rightarrow$ the torus is decomposed into periodic orbits with rational rotation number). Fix, once for all, one of these values.

\begin{figure}
\begin{center}
\includegraphics[width=12cm]{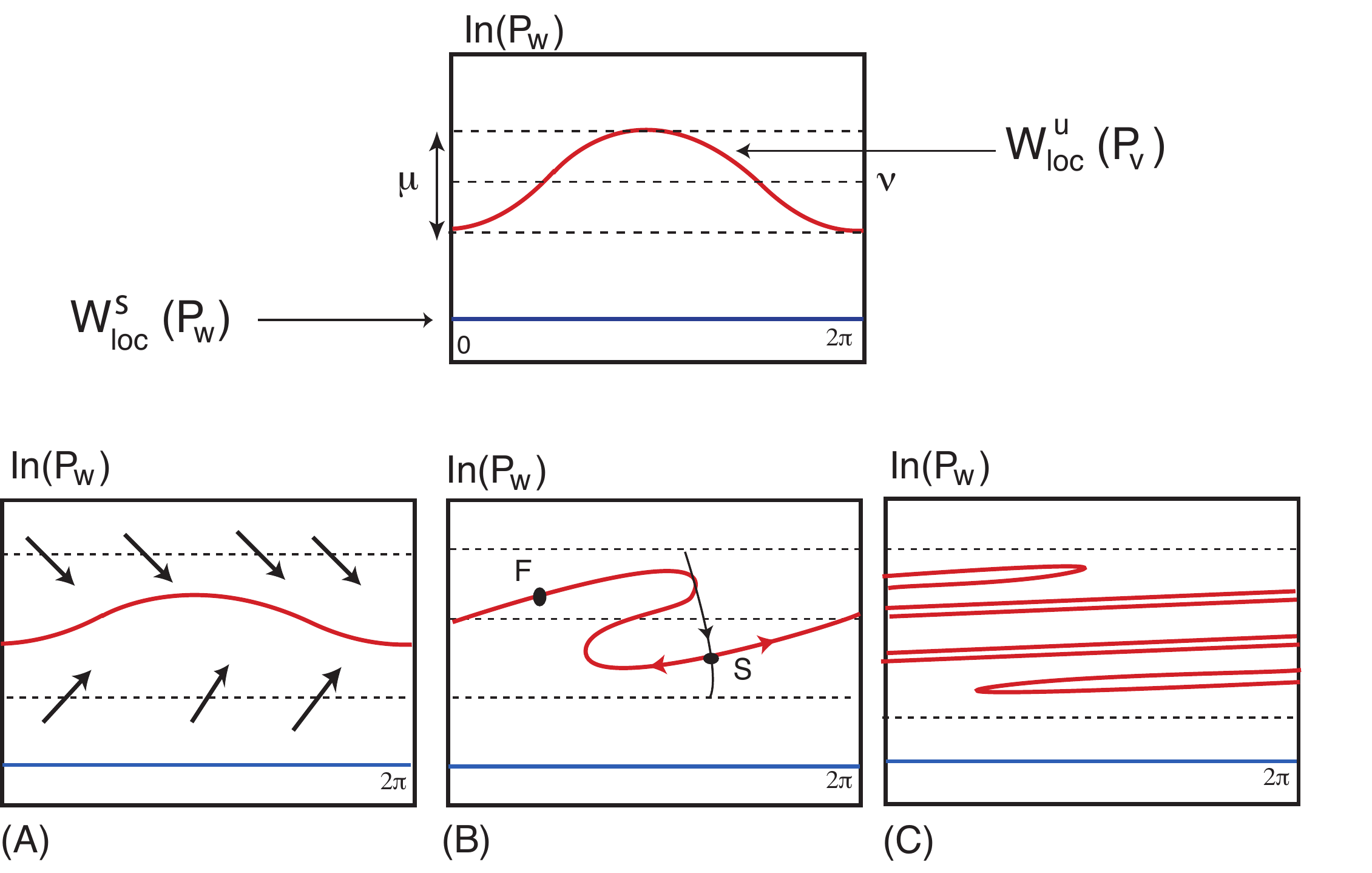}
\end{center}
\caption{\small  Image of $ \Rc(\mathcal{C})$ for different values of $\omega$ with  $\nu>0$ and $\mu$ fixed.
Transition from an invariant and attracting curve (A) to a rotational horseshoe (C), passing through a homoclinic tangency (B). One observes the breaking of the wave which accompanies the break of the invariant circle  (corresponding to the torus).  Here the point F is a sink and the point S is a saddle. In (C), a neighborhood of $r = 1$ is folded and mapped into itself, leading to the formation of rotational horseshoes. }
\label{Strange_attractors1_LR_torus}
\end{figure}

For such a $\nu_i>0$, increasing $\mu>0$ the coexistence of this pair of periodic orbits persists along a wedge, the so called \emph{Arnold tongue} \cite{Anishchenko}. As illustrated in Figure \ref{Strange_attractors1_LR_torus}, for a fixed $\mu>0$, we know that:
\begin{itemize}
\item $\Rc(\mathcal{C})$ is a closed curve on $Out^+(\Pp_\vv)$ because $\Rc|_\mathcal{D}$ is a diffeomorphism;
\item for $\omega \approx 0$,  this curve may be seen as the graph on $Out^+(\Pp_\vv)$ of a non-constant map  defined on $[0, 2\pi]$ (cf. \cite{Rodrigues2019}); the curve $\mathcal{C}$ is the $\omega$-limit of  $W^u(\Pp_\ww)$;
\item  by Lemma \ref{lema:spiral}, there exists $\omega_0>0$ and a segment $S\subset \mathcal{C}$ such that $\Rc (S)\cap \mathcal{D}$ is a curve containing a segment across $\mathcal{D}$.
\end{itemize}

This means that the  curve  $\mathcal{C}$ starts to develop folds as in Figure~\ref{Strange_attractors1_LR_torus} (B) and (C). 
If $\omega>\omega_0$  it creates the rotational horseshoes proved in \S \ref{prova_parte1}.
 Within this wedge, 
Anishchenko, Safonova and Chua have shown in \cite{Anishchenko}
 that there are curves on the parameter space $(\nu, \mu)$ corresponding to a quadratic homoclinic tangency associated to a dissipative periodic point of the first return map $\Rc$.   

Using now the results of Mora and Viana in  \cite{MV93}, there exists a positive measure set $\mathcal{U}$ of parameter values, so that for every $(\nu, \mu) \in \mathcal{U}$, the return map $\mathcal{R}$ 
admits a strange attractor  of H\'enon-type. 
These strange attractors are supported in SRB measures. 

\begin{remark}
In this type of result, the number of connected components with which the strange attractors intersect the section $\Sigma$ is not specified nor is the size of their basins of attraction. The strange attractors coexist with sinks from Newhouse phenomena. A discussion of these results may be found in \cite{CR2021, Rodrigues2019, TS86}. 
\end{remark}

\section{Proof of Proposition \ref{periodic_solution_prop} }
\label{PropA}

This proposition concerns the case $\mu=0$ when the time-periodic perturbation to $\dot x=\Fc_\nu(x)$ is constant.
For $\nu=0$, Properties~\ref{A1} to \ref{A4} and \ref{A6} of $\dot x=\Fc_0(x)$ were established in  \cite[Theorem 7]{ACL06}, with $\MM^2=\EU^2$.
 
\begin{figure}
\begin{center}
\includegraphics[width=65mm]{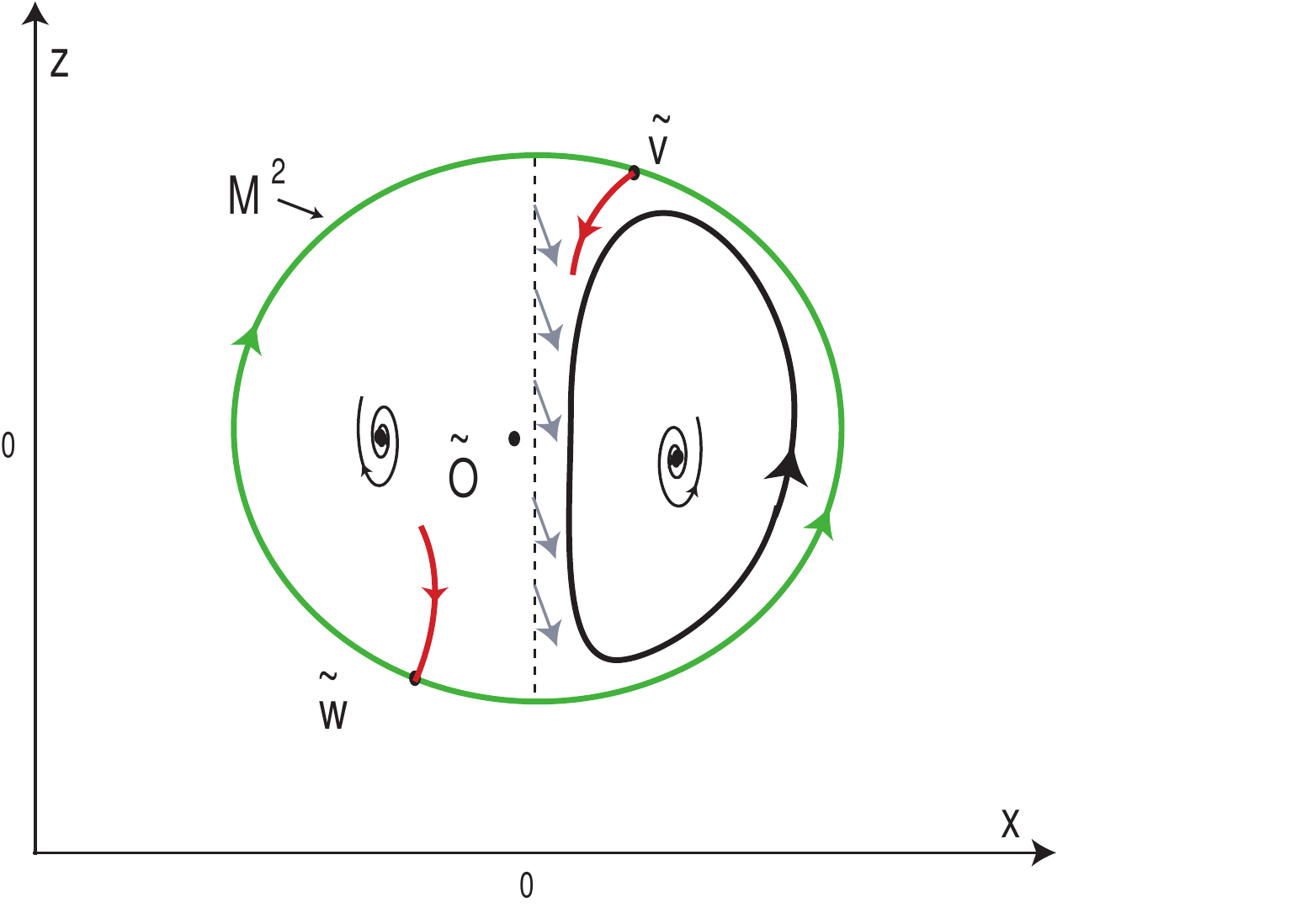}\includegraphics[width=65mm]{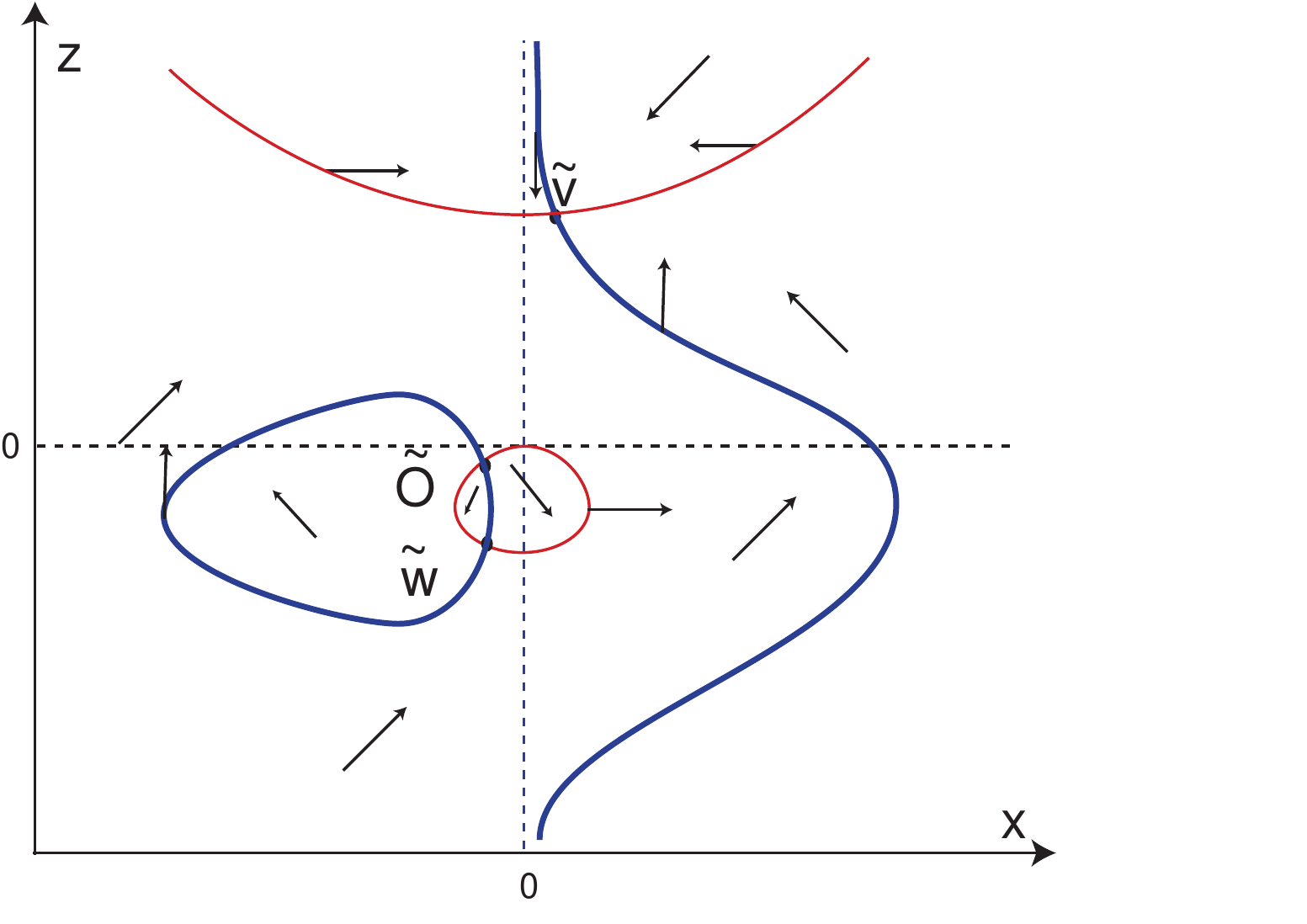} \\\includegraphics[width=65mm]{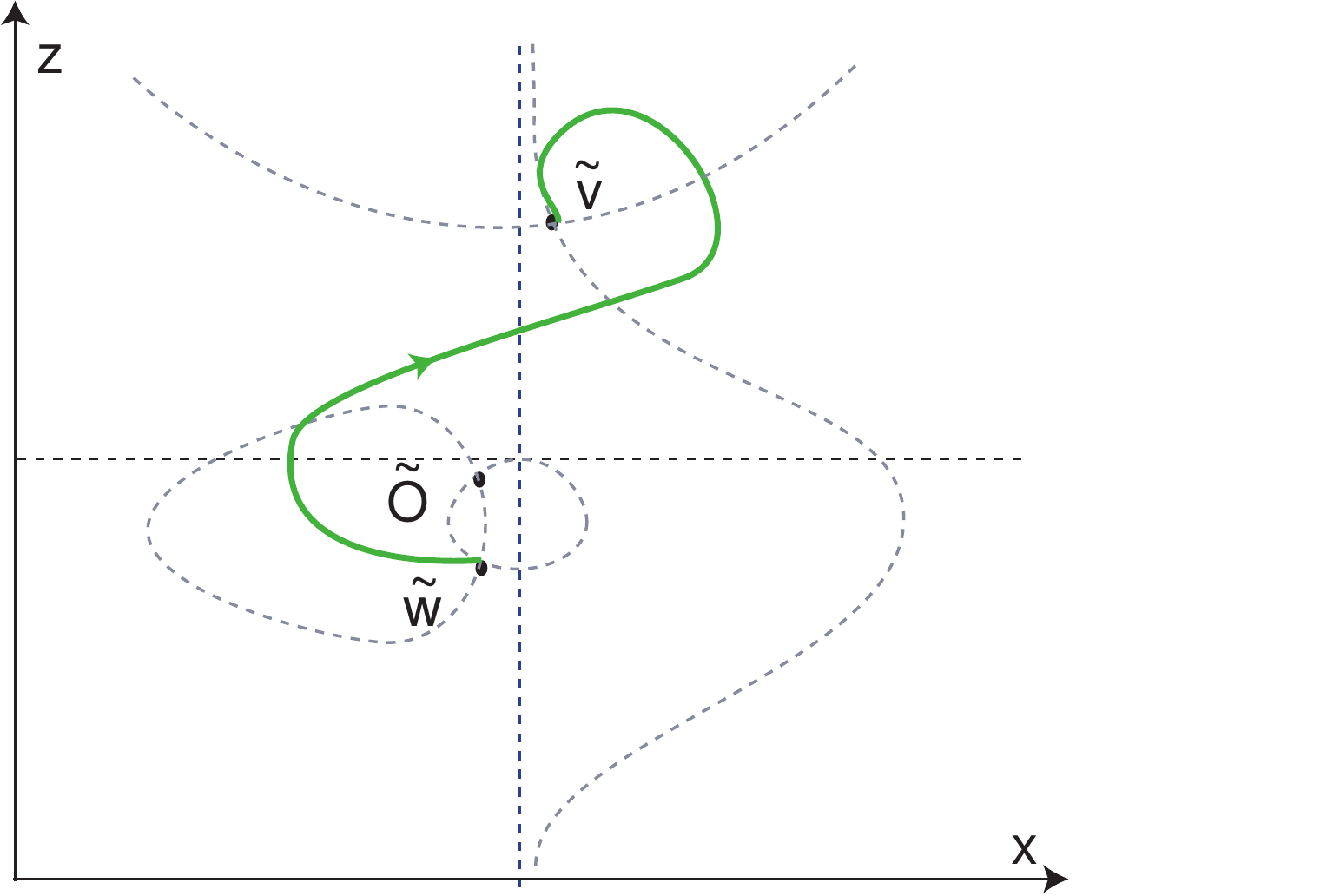}
\end{center}
\caption{\small Qualitative phase portrait for the dynamics or \eqref{general4} in $\MM^2$ with $y>0$, projected   into the $(x_1,0,x_3)$ plane for  $\mu=0$ and $\nu>0$ . }
\label{G1G2_5A}
\end{figure}

\begin{proof} The $\kappa$-equivariance is easily checked directly from the expression of $\Fc_\nu$.

 The first part of the proof of Proposition \ref{periodic_solution_prop}  consists in establishing that Properties~\ref{A1} to \ref{A3}  persist when the perturbation term $\nu(1-x_1)$ is added.
 Properties~\ref{A4} and \ref{A6} are established in \cite{ACL06} and Property~\ref{A5} is a consequence of their results.
Then it remains to show that two periodic solutions are created by the perturbation when the connections from $\ww$ to $\vv$ are broken, as stated in \ref{A7}.
Addressing the  persistence and property  \ref{A7} constitutes the remainder of this proof.

\begin{itemize}
\item[\ref{A1}] 
Since for $\nu=0$ the sphere $\EU^2$ is normally hyperbolic as in  \cite{HPS}, then for small $\nu\ne 0$ it persists as a flow-invariant, normally hyperbolic, globally attracting manifold $\MM^2$. See also the analysis by \cite{HG97}.  

\item[\ref{A2}] 
For  $\nu=0$ the only equilibria in the flow-invariant plane $\Fix(\ZZ_2(\kappa))$ are $\vv$ and $\ww$ above, as well as  the origin $O$.
Since they are hyperbolic, then their hyperbolic continuations should exist within the plane $\Fix(\ZZ_2(\kappa))$.
Another way to see this is to  solve
$$
\left\{\begin{array}{l}
F_1(x_1,x_3)= x_1[(1-x_1^2-x_3^2)-\alpha  x_3 +\beta x_3^2 - \nu] + \nu =0\\ \\
F_2(x_1,x_3)=  x_3[(-x_1^2-x_3^2) -\beta x_1^2] +\alpha x_1^2=0 .
\end{array}\right.
$$

\begin{figure}
\begin{center}
\includegraphics[width=75mm]{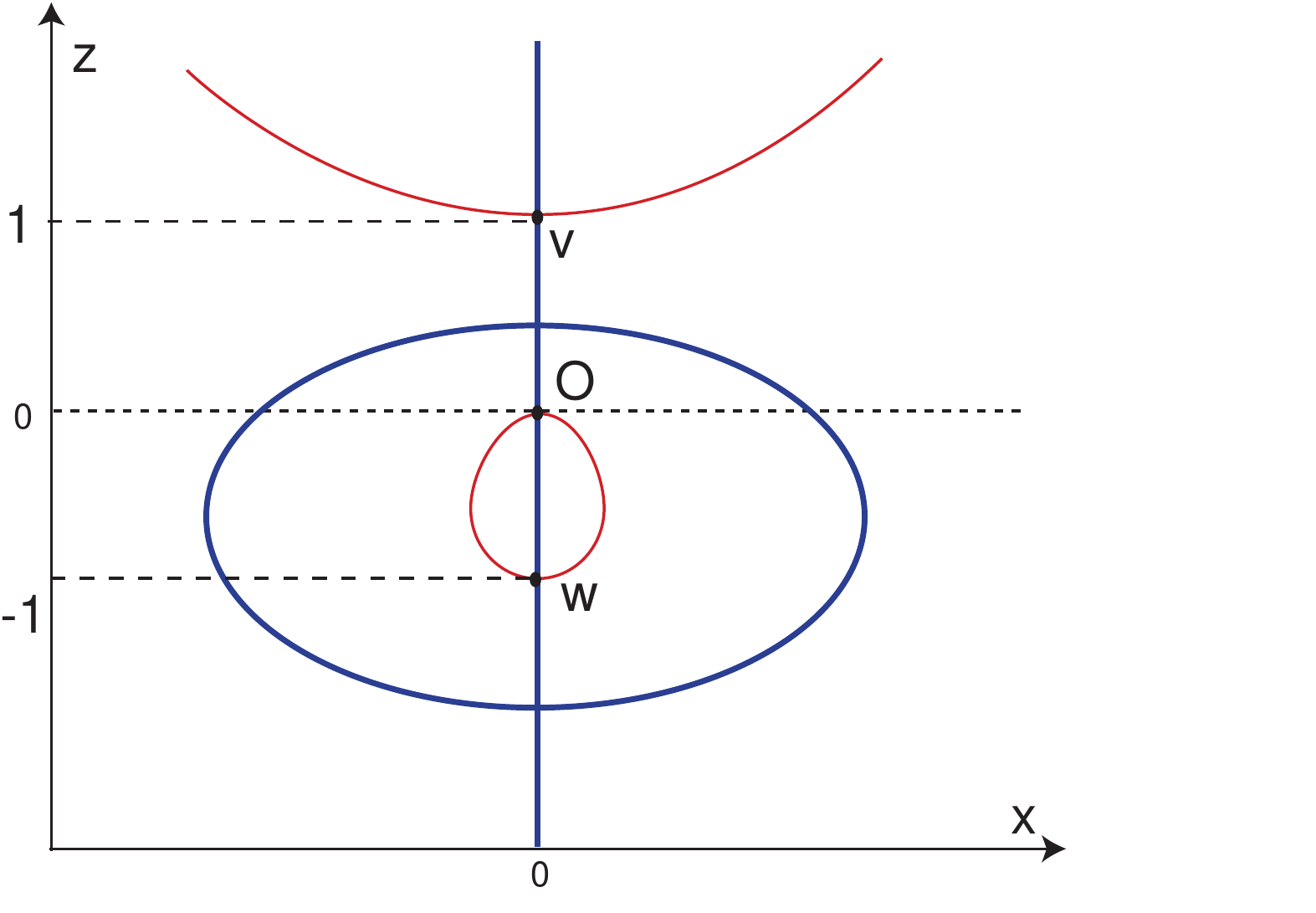}\qquad \includegraphics[width=75mm]{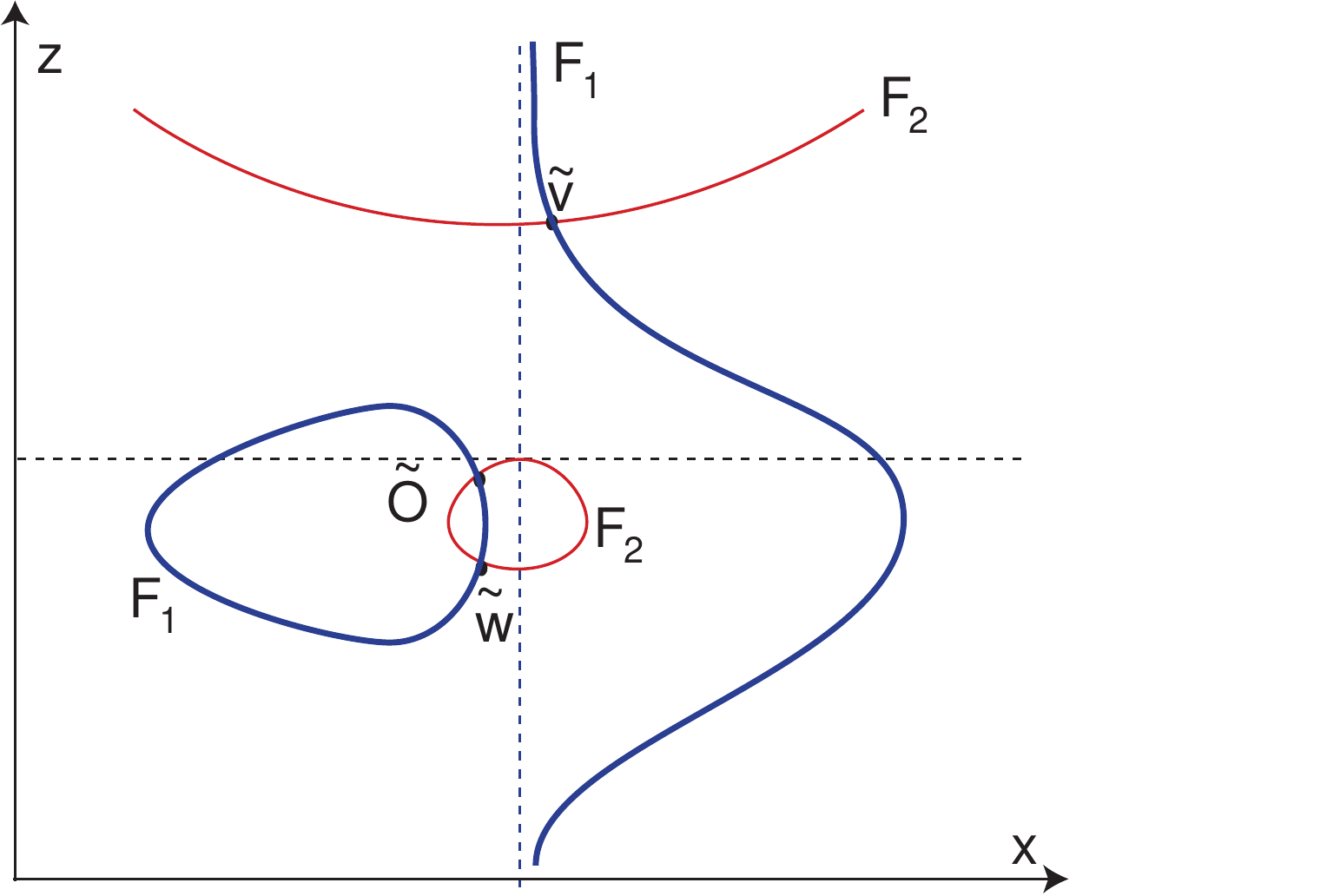}
\end{center}
\caption{\small Equilibria of $\dot\zeta=\Fc_{\nu}(\zeta)$ for $\nu=0$ (left) and $\nu=0.5$ (right), $\alpha=1$, $\beta=-0.1$ in the flow-invariant subspace $\Fix(\ZZ_2(\kappa))$ occur at the intersection of the curves $F_1(x_1,x_3)=0$ (blue) and $F_2(x_1,x_3)=0$ (red), here plotted with Maxima.
}
\label{G1G2}
\end{figure}

Figure~\ref{G1G2} shows the curves $F_1=0$ and $F_2=0$ plotted with Maxima. 
For $\nu=0$ the curves intersect transversely at $O$, $\vv$ and $\ww$ and this property is preserved for small $\nu\ne 0$.
\item[\ref{A3}] 
Within the flow-invariant plane $\Fix(\ZZ_2(\kappa))$, the origin is a repelling source,  $\vv$ is a sink and $\ww$  is a saddle. Since both the plane $y=0$ and the manifold $\MM^2$ are flow-invariant, this means that there are two heteroclinic connections from 
$\tilde{\ww}$ to $\tilde{\vv}$.

\item[\ref{A5}] 
In \cite{ACL06} it is established that the only other equilibria in $\EU^2$ are the four hyperbolic repelling foci $(\pm \sqrt{2}/2, \pm \sqrt{2}/2, 0)$. 
This means that for  $\nu=0$, by the Poincar\'e-Bendixon Theorem,  the $\omega$-limit set of all other points in $\EU^2$ must be contained in the heteroclinic cycles that contain $\vv$ and $\ww$.
The unstable foci  remain for $\nu\ne 0$ small. 

\item[\ref{A7}]
The flow-invariant subspace $\Fix(\ZZ_2(\kappa))$ divides $\MM^2$ in two flow-invariant components.
We will show that the $x_2>0$ component contains a non-constant periodic solution, the proof for $x_2<0$ follows from the symmetry.
For $x_1=0$ and $\mu=0$  the expression \eqref{general4} yields $\dot x=\nu\ne 0$.
Suppose $\nu>0$, then the region $x_1>0$, $x_2>0$ in $\MM^2$ is positively invariant, see Figure~\ref{G1G2_5A}.
This region only contains one equilibrium, one of the repelling foci in \ref{A4}, hence by the Poincar\'e-Bendixon Theorem,
the $\omega$-limit of the unstable manifold of $\tilde\vv$ is an attracting periodic solution. 
When $\nu<0$  the periodic solution appears  for $x_1<0$.   
The period tends to $+\infty$ as $\nu$ goes to $0$, since the periodic trajectory  accumulates on
 the heteroclinic cycle.
\end{itemize}
\end{proof}

Proposition \ref{periodic_solution_prop} shows that for sufficienly small $|\nu|>0$, each heteroclinic cycle that occurred in the fully symmetric case  
is replaced by a stable hyperbolic periodic solution. Using the reflection symmetry $\ZZ_2(\kappa)$, two stable periodic solutions co-exist, one in each connected component of $\mathcal{M}\backslash \Fix(\ZZ_2(\kappa))$. Their period tends to $\infty$ when $\nu$ vanishes  and their basin of attraction must contain the basin of attraction of $\Sigma_0$.

\section{Discussion and concluding remarks}
\label{discussion}

    \emph{Routes  to chaos}  have been a recurrent concern on nonlinear dynamics during the last decades \cite{AS91}. The novelty of the present paper is the illustration of a \emph{new route} for the emergence of strange attractors  from an attracting heteroclinic network as a codimension-two phenomenon. 

\subsection{Literature}
In \cite{Kaneko}  Kaneko investigates numerically the
 bifurcations of tori in a two-parameter family of dissipative coupled maps. Each map undergoes a period-doubling cascade accumulating on a given parameter value, which generates chaos in the coupled system.   In the same setting Bakri and Verhulst show in  \cite{Bakri}  that,  for small amplitudes, zero damping and zero coupling, their system reveals a periodic solution which undergoes a Hopf bifurcation, generating an attracting torus. 
 They use numerical bifurcation techniques to show   how the torus gets destroyed by dynamical and topological changes in the involved manifolds. The results agree with \cite{Ruelle, TS86}.

In the context of dissipative vortex  dynamics Fleurantin and James have studied in \cite{FJ2020} the  \emph{Langford system}, a   one-parameter family of three-dimensional vector fields.
The flow of this model exhibits a sink, two saddle-foci of different Morse indices and a non-trivial periodic solution with a complex conjugate pair of Floquet exponents. The frequency of the periodic solution together with the frequency of the complex exponent constitute two competing natural modes of oscillation. 
The periodic orbit undergoes a 
bifurcation  giving rise to observable chaos  through the same mechanism of \cite{Bakri}.
 These authors studied the evolution of the torus,  its loss of differentiability,  and the appearance of a strange attractor via the existence of  tangencies.  
  The relative position of the manifolds according to the parameter allows the authors to prove the existence of \emph{bistability} between an equilibrium and a torus.  

\subsection{Heteroclinic tangle}
The formation of the horseshoe of Theorem \ref{role_omega} has a different nature to   those found in \cite{ACL06, LR17} -- in this case, the shift dynamics is obtained via the transverse intersection of two two-dimensional invariant manifolds.  The parameter $\omega$ is not necessary to prove the existence of chaos.  The non-wandering set associated to the network contains, but does not coincide with, the suspension of horseshoes; it contains  infinitely many heteroclinic pulses and attracting limit cycles with long  periods, coexisting with sets with positive entropy, giving rise to the so called \emph{quasi-stochastic attractors} \cite{LR17}.  The sinks  have long periods and narrow basins of attraction.

\subsection{Open questions}
In the context of this class of examples, some problems remain to be solved:
\begin{enumerate}
\item the basins of attraction of the strange attractors of Theorem B are relatively small in terms of Lebesgue measure (they are close to Newhouse domains). Could we improve Theorem B in order to get the existence of a ``larger'' strange attractor?
\item is it possible to generalize our result for clean heteroclinic networks (networks whose unstable manifolds are contained within it) whose connections are one-dimensional?
\end{enumerate}
We believe that these problems can be tackled by using the theory of rank-one attractors developed by Wang and Young \cite{WY}. We defer these tasks for future work.


\begin{thebibliography}{99}
\bibitem{AHL2001} V.S. Afraimovich, S-B Hsu, H. E. Lin,\emph{Chaotic behavior of three competing species of May-Leonard model under small periodic perturbations}. Int. J. Bif. Chaos, 11(2), 435--447, 2001
\bibitem{AS91} V. S. Afraimovich,  L. P.  Shilnikov, \emph{On invariant two-dimensional tori, their breakdown and stochasticity, Methods of the Qualitative Theory of Differential Equations } (Gor'kov. Gos. University) [1983], 3--26. Translated in: Amer. Math. Soc. Transl., 149:2, 201--212, 1991
%
\bibitem{Anishchenko} V. Anishchenko, M. Safonova,  L.  Chua,  \emph{Confirmation of the Afraimovich-Shilnikov torus-breakdown theorem via a torus circuit}, IEEE Transactions on Circuits and Systems I: Fundamental Theory and Applications 40:11, 792--800, 1993
%
%
\bibitem{ACL06} M.A.D. Aguiar, S.B.S.D. Castro,  I. S. Labouriau, \emph{Simple Vector Fields with Complex Behavior}, {Int. J.  Bif. and Chaos}, {Vol. \textbf{16}} {No. \textbf{2}},  369--381, 2006
%
\bibitem{Bakri} T. Bakri, F. Verhulst,  \emph{Bifurcations of quasi-periodic dynamics: torus breakdown}, Zeitschrift fur angewandte Mathematik und Physik 65.6  1053--1076, 2014
%
\bibitem{CR2021} M. L. Castro, A. A. P. Rodrigues, \emph{Torus-breakdown near a heteroclinic attractor: a case study}, Int. J. Bif. Chaos, Vol. 31, No. 10 (2021) 2130029 (20 pages), DOI: 10.1142/S0218127421300299
%
\bibitem{BIRR} P. G. Barrientos, S. Ib\'a\~nez, A. A. Rodrigues, J. A. Rodr\'iguez, \emph{Emergence of Chaotic Dynamics from Singularities}, 32th Brazilian Mathematics Colloquium,  Instituto Nacional de Matem\'atica Pura e Aplicada (IMPA), Rio de Janeiro, 2019. xi+200 pp. ISBN: 978-85-244-0430-6, 2019
%
\bibitem{DT3} J. Dawes, T.-L. Tsai, \emph{Frequency locking and complex dynamics near a periodically forced robust heteroclinic cycle},  Phys. Rev. E, 74 (055201(R)), 2006
%
\bibitem{FJ2020} E. Fleurantin,  J. M. James, \emph{Resonant tori, transport barriers, and chaos in a vector field with a Neimark-Sacker bifurcation}. Communications in Nonlinear Science and Numerical Simulation, 85, 105226, 2020
%
\bibitem{GH} J. Guckenheimer,  P. Holmes, \emph{Nonlinear Oscillations, Dynamical Systems, and Bifurcations of Vector Fields, Applied Mathematical Sciences} (Springer-Verlag), 42, 1983
%
\bibitem{Herman} M. Herman,  \emph{Mesure de Lebesgue et Nombre de Rotation, Lecture Notes in Math} (Springer), 597, 271--293, 1977
%
\bibitem{HPS} M. Hirsch, C. Pugh, M. Shub, \emph{Invariant Manifolds}, Springer-Verlag, 1977
\bibitem{HG97} C. Hou, M. Golubitsky, \emph{An example of symmetry breaking to heteroclinic cycles}, J. Diff. Eqs., 133:30--48, 1997
\bibitem{Kaneko} K. Kaneko, \emph{Doubling of torus},  Progress of theoretical physics 69.6: 1806--1810, 1983
\bibitem{KM1} M. Krupa, I. Melbourne, \emph{Asymptotic Stability of Heteroclinic Cycles in Systems with Symmetry, }Ergodic Theory and Dynam. Sys., Vol. {15}, 121--147, 1995
\bibitem{KM2} M. Krupa, I. Melbourne, \emph{Asymptotic Stability of Heteroclinic Cycles in Systems with Symmetry,\ II, }Proc. Roy. Soc. Edinburgh, 134A, {1177--1197},  2004
\bibitem{LR17}  I. S. Labouriau, A. A. P. Rodrigues, \emph{On Takens' last problem: Tangencies and time averages near heteroclinic networks}, \emph{Nonlinearity}, \textbf{30}, 1876--1910, 2017
 \bibitem{LR18}    I.S. Labouriau, A.A.P. Rodrigues, \emph{Bifurcations from an attracting heteroclinic cycle under periodic forcing}, J. Differential Equations, \textbf{269}, 4137--4174, 2020
%
\bibitem{MV93} L. Mora, M.  Viana,  \emph{Abundance of strange attractors, Acta Math.} 171(1)  1--71, 1993
%
\bibitem{Passeggi} A. Passeggi, R. Potrie,  M.  Sambarino, \emph{Rotation intervals and entropy on attracting annular continua, Geometry \& Topology} 22:4, 2145--2186, 2018
%
\bibitem{Rabinovich06} M.I. Rabinovich, R. Huerta,  P. Varona, \emph{Heteroclinic synchronization: Ultra-subharmonic locking}, Phys. Rev. Lett., 96:014101, 2006
%
\bibitem{Rodrigues2019} A. A. P. Rodrigues,  \emph{Unfolding a Bykov attractor: from an attracting torus to strange attractors, Journal of Dynamics and Differential Equations}  2020  https://doi.org/10.1007/s10884-020-09858-z
%
%
\bibitem{Shilnikov_book_1} L. P. Shilnikov, A. Shilnikov, D. Turaev, L. Chua. \emph{ Methods Of Qualitative Theory In Nonlinear Dynamics} (Part I). World Sci. Singapore, New Jersey, London, Hong Kong, 1998.
%
\bibitem{Ruelle} D. Ruelle,   \emph{Differentiable dynamical systems and the problem of turbulence}, Bulletin of the American Mathematical Society 5.1,  29--42, 1981
\bibitem {TD1} T.-L. Tsai, J. Dawes, \emph{Dynamic near a periodically-perturbed robust heteroclinic cycle}, Physica D, 262, 14--34, 2013
%
\bibitem{TS86} D. V. Turaev,  L. P. Shilnikov,  \emph{Bifurcation of torus-chaos quasi-attractors, Mathematical Mechanisms of Turbulence (Russian)}   Akad. Nauk Ukrain. SSR, Inst. Mat., Kiev: 113--121, 1986.
%
\bibitem{Wang_2013} Q. Wang, \emph{Heteroclinic tangles in time-periodic equations}, J. Diff. Eqs. 254  1137--1171, 2013  
%
\bibitem{WY} Q. Wang, L.S.  Young.  \emph{From Invariant Curves to Strange Attractors}, Commun. Math. Phys.  225--275, 2002.
%
%
\bibitem {Wiggins} S.Wiggins, \emph{Introduction to applied nonlinear dynamical systems and chaos}, Springer-Verlag, 2003
\end{thebibliography}
\end{document}